\documentclass[a4paper,12pt]{article}

\usepackage{amssymb,amsmath,amsthm}
\usepackage[totalwidth=17cm,totalheight=24cm]{geometry}
\usepackage{graphicx}

\newcommand{\C}{{\mathbb C}}
\newcommand{\R}{{\mathbb R}}
\newcommand{\Z}{{\mathbb Z}}
\newcommand{\Hq}{{\mathbb H}}
\newcommand{\Oc}{{\mathbb O}}

\newcommand{\id}{{\mathbb I}}
\newcommand{\im}{{\rm i\,}}

\newcommand{\e}{{\bf e}}

 \theoremstyle{plain}
  \newtheorem{theorem}{Theorem}
   \newtheorem{proposition}{Proposition}

  \theoremstyle{definition}
  \newtheorem{definition}[theorem]{Definition}

  \theoremstyle{remark}

\newcommand{\be}{\begin{eqnarray}}
\newcommand{\ee}{\end{eqnarray}}

\begin{document}
 \pagestyle{plain}
\title{Lorentzian Cayley Form}

\author{Kirill Krasnov\\ {}\\
{\it School of Mathematical Sciences, University of Nottingham, NG7 2RD, UK}}

\date{April 2023}
\maketitle

\begin{abstract}\noindent Cayley 4-form $\Phi$ on an 8-dimensional manifold $M$ is a real differential form of a special algebraic type, which determines a Riemannian metric on $M$ as well as a unit real Weyl spinor. It defines a ${\rm Spin}(7)$ structure on $M$, and this ${\rm Spin}(7)$ structure is integrable if and only if $\Phi$ is closed. We introduce the notion of a complex Cayley form. This is a one-parameter family of complex 4-forms $\Phi_\tau$ on $M$ of a special algebraic type. Each $\Phi_\tau$ determines a {\it real} Riemannian metric on $M$, as well as a complex unit Weyl spinor $\psi_\tau$. The subgroup of ${\rm GL}(8,\R)$ that stabilises $\Phi_\tau, \tau\not=0$ is ${\rm SU}(4)$, and $\Phi_\tau$ defines on $M$ an ${\rm SU}(4)$ structure. We show that this ${\rm SU}(4)$ structure is integrable if and only if $\Phi_\tau$ is closed. 

We carry out a similar construction for the split signature case. There are now two one-parameter families of complex Cayley forms. A complex Cayley form of one type defines an ${\rm SU}(2,2)$ structure, a form of the other type defines an ${\rm SL}(4,\R)$ structure on $M$. As in the Riemannian case, these structures are integrable if and only of the corresponding complex Cayley forms are closed. Our central observation is that there exists a special member of the second one-parameter family of complex Cayley forms, which we call the Lorentzian Cayley form. This 4-form has the property that it is calibrated by Lorentzian 4-dimensional subspaces $H,H^\perp$. In particular, in a basis adapted to such a calibration, the Lorentzian Cayley form is built from the complex self-dual 2-forms for $H,H^\perp$. We explain how these observations solve a certain puzzle that existed in the context of 4-dimensional Lorentzian geometry. 
\end{abstract}

\section{Introduction}

A real 4-form $\Phi\in \Lambda^4 M$ on an 8-dimensional manifold $M$ is called {\bf Cayley} form if it admits a compatible with it Riemannian metric $g_\Phi \in S^2 T^*M$, see Section \ref{sec:definitions} for the precise definition. The metric $g_\Phi$ is determined uniquely by $\Phi$. Given the metric $g_\Phi$, there is the associated Clifford algebra $\Gamma: TM\to {\rm End}(S)$, where $S=S^+\oplus S^-$ is the space of spinors, and the Clifford matrices satisfy their usual defining relations $\Gamma(\xi)\Gamma(\xi) = g_\Phi(\xi,\xi) \id, \xi\in TM$. Further, given $\Phi$ and $g_\Phi$ there exists a special real unit Weyl spinor $\psi_\Phi \in S^+$ such that 
\be\label{Phi-spinor}
\langle \psi_\Phi, \Gamma(\xi_1) \Gamma(\xi_2) \Gamma(\xi_3) \Gamma(\xi_4) \psi_\Phi \rangle = \Phi(\xi_1,\xi_2,\xi_3,\xi_4).
\ee
Here $\xi_{1,2,3,4}\in TM$ and $\langle\cdot,\cdot\rangle$ is the ${\rm Spin}(8)$ invariant inner product on $S^+, S^-$. In other words, a Cayley form $\Phi$ determines both the metric $g_\Phi$ as well as a unit real Weyl spinor $\psi_\Phi$ such that the Cayley form itself can be reproduced as an appropriate square (\ref{Phi-spinor}) of this spinor.

This statement can be confirmed by a dimensions count. A Cayley form is a 4-form on $\Lambda^4 M$ of a special algebraic type, whose stabiliser in ${\rm GL}(8,\R)$ is ${\rm Spin}(7)$. Note that ${\rm Spin}(7)$ is also the stabiliser in ${\rm Spin}(8)$ of a real unit Weyl spinor. The space of Cayley forms is then 
\be\label{dim-Cayley}
\Lambda^4(\R^8) \supset \{ \text{Cayley forms}\} = {\rm GL}(8,\R)/{\rm Spin}(7).
\ee
This is the space of dimension $64-21=43$. This is also the dimension $36$ of the space of metrics on $\R^8$ plus the dimension $7$ of the space of real unit spinors, so we can write
\be
\{ \text{Cayley forms}\}  = \{ \text{metrics on} \,\, \R^8\} \oplus \{ \text{unit Weyl spinors}\}.
\ee

It is useful to write the formula for $g_\Phi$ explicitly. When $\Phi$ is a 4-form of the Cayley algebraic type we have 
\be\label{metric-Phi}
\frac{1}{6} \xi_1\lrcorner \xi_2\lrcorner \Phi \wedge \eta_1\lrcorner \eta_2\lrcorner \Phi \wedge \Phi = (g_\Phi(\xi_1,\eta_1)g_\Phi(\xi_2,\eta_2) - g_\Phi(\xi_1,\eta_2)g_\Phi(\xi_2,\eta_1)) v_\Phi. 
\ee
where $\xi_{1,2}, \eta_{1,2} \in TM$ and $v_\Phi$ is the volume form of the metric $g_\Phi$ determined by $\Phi$. This is one of the ways in which the metric $g_\Phi$ is characterised by $\Phi$.  

A Cayley form $\Phi$ can be written in a way that exhibits calibrations. We remind that a 4-dimensional subspace $H\subset TM$ is called Cayley, or {\bf calibrated} by $\Phi$, if for any orthonormal basis $\xi_{1,2,3,4}\in H$ (with respect to $g_\Phi)$ we have $\Phi(\xi_1,\xi_2,\xi_3,\xi_4)=\pm 1$. The orthogonal complement $H^\perp$ of a Cayley subspace is also Cayley, and $\Phi$ can be written in a way that exhibits the ${\rm SU}(2)\times{\rm SU}(2)\times{\rm SU}(2)$ subgroup of its stabiliser ${\rm Spin}(7)$. To this end one introduces an orthonormal basis $\Sigma^i,\Sigma^{'i}, i=1,2,3$ in the space of self-dual 2-forms in $H,H^\perp$. The Cayley form then takes the form (\ref{Phi-calibr}).

Another way to think about this geometric setup is that a Cayley form $\Phi$ on $M$ (together with a choice of a unit vector) identifies the tangent space $T_p M$ at each point with the space of octonions $\Oc$ (so that the chosen vector is identified with the identity octonion). Indeed, as we shall review in Section \ref{sec:definitions}, a Cayley form equips $T_p M$ with what is known as the triple cross product. This triple cross product comes from the octonion product once the tangent space $T_p M$ is identified with $\Oc$, see below. From this octonionic point of view, a calibration is a choice of a copy of quaternions $\Hq\subset \Oc$. 

The above discussion assumed that the metric $g_\Phi$ determined by $\Phi$ is Riemannian. However, a real 4-form of the algebraic type that admits a compatible metric can give rise to both a Riemannian and a split signature metrics via (\ref{metric-Phi}). In the latter case $\Phi$ is a real 4-form of a special algebraic type whose ${\rm GL}(8,\R)$ stabiliser is ${\rm Spin}(4,3)$, and $\Phi$ determines a split signature metric $g_\Phi$ as well as a real unit Weyl spinor $\psi_\Phi$ of ${\rm Spin}(4,4)$ such that (\ref{Phi-spinor}) continues to hold. A Cayley form of split type is a real 4-form that identifies $T_p M$ with the space of split octonions. 

The first new observation of this paper is that we can extend the notion of 4-forms of Cayley type also to {\bf complex-valued} 4-forms on $\Lambda^4 T^*M$. We will call a complex 4-form $\Phi_\tau\in \Lambda^4_\C T^*M$ a {\bf complex Cayley} form, or a complex form of a Cayley type, if there is a compatible with it, in the sense of formula (\ref{metric-Phi}), {\bf real} metric. Examples of such complex 4-forms are easy to produce by taking the real metric of signature either all plus or split, building the corresponding Clifford algebra, and then applying the formula (\ref{Phi-spinor}) with a {\bf complex} unit spinor $\psi$. The resulting complex-valued 4-form then continues to produce a real metric via (\ref{metric-Phi}).

Interestingly, dropping the reality requirement for $\psi_\Phi$ in (\ref{Phi-spinor}) allows to construct complex Cayley forms in eight dimensions in any signature. However, the most interesting cases to consider are still Riemannian and split, where real spinors exist and complex Cayley forms can be compared to their real cousins. In the case of Riemannian signature, there is a single orbit in $S^+$ of complex unit spinors $\psi_\tau: \langle \psi_\tau, \psi_\tau \rangle =1$ for each value of $\langle \hat{\psi}_\tau,\psi_\tau\rangle:=\cosh(2\tau)$, where the hat in $\hat{\psi}_\tau\in S^+$ is the ${\rm Spin}(8)$ invariant complex conjugation. The ${\rm Spin}(8)$ stabiliser of a complex unit spinor $\psi_\tau, \tau\not=0$ is ${\rm SU}(4)$. This will also be the ${\rm GL}(8,\R)$ stabiliser of the complex Cayley form $\Phi_\tau$ produced from $\psi_\tau$. Thus, a complex Cayley form that gives rise to a Riemannian signature metric $g_{\Phi_\tau}$ endows $M$ with a ${\rm SU}(4)$ structure rather than ${\rm Spin}(7)$ structure produced on $M$ by a real Cayley form. It is illuminating to write an explicit formula for this complex Cayley form. It is given by
\be\label{Phi-c}
\Phi_\tau= \langle \psi_\tau, \Gamma\Gamma\Gamma\Gamma \psi_\tau \rangle =\cosh(2\tau){\rm Re}(\Omega)-\frac{1}{2} \omega\wedge \omega + \im \sinh(2\tau){\rm Im}(\Omega).
\ee
Here $\Omega$ is the $(4,0)$ form for the complex structure $J$ that the ${\rm SU}(4)$ structure $\Phi_\tau$ defines, and $\omega$ is the $(1,1)$ Kahler form. This should be compared with the expression (\ref{Phi-compl-form}) for the usual real Cayley form, once an ${\rm SU}(4)$ structure is chosen. Both (\ref{Phi-c}) and  (\ref{Phi-compl-form}) are manifestly ${\rm SU}(4)$ invariant. When $\tau=0$ the complex Cayley form becomes real and turns into the usual real Cayley form 
$\Phi_0=\Phi$. The invariance of $\Phi_\tau$ at $\tau=0$ is thus increased from ${\rm SU}(4)$ to ${\rm Spin}(7)$. 

There is the dimensions count similar to (\ref{dim-Cayley}) but for complex Cayley forms. 
\be\label{dim-Cayley-compl}
\Lambda^4_\C(\R^8) \supset \{ \text{complex Cayley forms}\} = {\rm GL}(8,\R)/{\rm SU}(4).
\ee
This is the space of dimension $64-15=49$. This is also the dimension $36$ of the space of metrics on $\R^8$ plus the dimension $13$ of the space of complex unit Weyl spinors of fixed norm squared, so we can write
\be\nonumber
\{ \text{complex Cayley forms}\}  = \{ \text{metrics on} \,\, \R^8\} \oplus \{ \text{complex} \,\,\psi: \langle \psi,\psi\rangle=1, \langle\psi,\hat{\psi}\rangle = \cosh(2\tau)\}.
\ee

Complex Cayley forms of the type (\ref{Phi-c}) are interesting because they can be used to encode ${\rm SU}(4)$ structures in 8 dimensions in a single object $\Phi_\tau$. We remind the reader that the ${\rm SU}(4)$ structure on $M$ is a reduction of the principal ${\rm GL}(8,\R)$ frame bundle to an ${\rm SU}(4)$ subbundle.  One obtains such a reduction when the complex Cayley form is closed $d\Phi_\tau=0$, and every Calabi-Yau structure on $M$ defines a closed complex Cayley form via (\ref{Phi-c}). Proposition \ref{prop:complex-integrable} below is a statement to this effect. 

Similar to real Cayley forms, there is also a way of writing the complex Cayley form in a way adapted to a calibration. To describe this we need to generalise the notion of the calibration slightly, and say that $H\subset TM$ is calibrated by $\Phi_\tau$ if for any orthonormal (with respect to the real metric defined by $\Phi_\tau$) basis $\xi_{1,2,3,4}\in H$ the quantity $\Phi_\tau(\xi_1,\xi_2,\xi_3,\xi_4)=\pm \cosh(2\tau)$. The complex Cayley form $\Phi_\tau$ can then be written in a way adapted to one of its calibrations. It is no longer possible to write an expression that involves just the self-dual 2-forms for $H,H^\perp$, because the stabiliser subgroup of $\Phi_\tau$ is smaller. The expression that arises is given by (\ref{Phi-tau-calibr}) and involves an orthonormal basis of 1-forms for both $H,H^\perp$. This way of writing $\Phi_\tau$ exhibits the ${\rm SO}(4)$ subgroup of its stabiliser. This is the diagonal group of rotations that acts simultaneously on $H,H^\perp$. 

The second and main new observation of this paper is that there also exist complex Cayley forms that are compatible with the split signature metric on $M$, and that in the split case the story is richer and exhibits phenomena that do not have analogs in the Riemannian case. The story in the split signature case starts analogously to that in the Riemannian case. One introduces complex Cayley forms $\Phi_\tau \in \Lambda^4_\C TM$ to be those that are compatible with a real split signature metric on $M$. They are again constructed from complex unit spinors $\psi\in S^+$ of ${\rm Spin}(4,4)$ via (\ref{Phi-spinor}). The new situation is that there are now several different types of complex unit spinors, with different stabiliser. In particular, there are two one-parameter families of orbits, distinguished by the value of $\langle \hat{\psi},\psi\rangle$. Analogously to the Riemannian case, there is a family $\Phi_\tau$, where $\langle \hat{\psi},\psi\rangle =\cosh(2\tau)$. The stabiliser of a generic 4-form $\Phi_\tau$ is ${\rm SU}(2,2)$, and so $\Phi_\tau$ defines a real split signature metric on $M$, as well as an orthogonal complex structure, together with the top holomorphic form $\Omega$. The special member of this family of 4-forms $\Phi_0$ coincides with the real split Cayley form, whose stabiliser is ${\rm Spin}(4,3)$. 

The other arising one-parameter family of complex Cayley forms $\Phi_\theta, \theta\in[0,\pi/2]$ does not have a Riemannian analog. It defines a real split signature metric as well as a complex unit Weyl spinor $\psi_\theta$ with $\langle \hat{\psi}_\theta,\psi_\theta\rangle=\cos(2\theta)$. A generic member of this family of 4-forms has the stabiliser ${\rm SL}(4,\R)$. This means that $\Phi_\theta$ endows $M$ with a para-complex structure (see below for the definition) whose real null four-dimensional subspaces are those on which the stabiliser ${\rm SL}(4,\R)$ acts. The 4-form $\Phi_\theta$ can be written in terms of the data that the para-complex structure defines. Let $\Omega_\pm$ be the real decomposable 4-forms given by the products of eigendirections of the para-complex structure (these are the real analogs of $\Omega,\bar{\Omega}$ in the case of complex structure), and $\omega_r$ be the analog of the Kahler form. We then have
\be
\Phi_\theta = \cos(2\theta) \frac{1}{2}(\Omega_++\Omega_-) + \im \sin(2\theta) \frac{1}{2}(\Omega_+ - \Omega_-) + \frac{1}{2} \omega_r \wedge \omega_r.
\ee
For $\theta=0,\pi/2$ the 4-form $\Phi_\theta$ is again real, and coincides with the real split Cayley form with stabiliser ${\rm Spin}(4,3)$. 

Our main observation is that there exists yet another special member of this second family of complex Cayley forms. It corresponds to $\theta=\pi/4$. The corresponding complex unit spinor is orthogonal to its complex conjugate $\langle \hat{\psi}_{\pi/4},\psi_{\pi/4}\rangle=0$. We will refer to complex Cayley form $\Phi_L:=\Phi_{\pi/4}$ as the {\bf Lorentzian Cayley form}. The explicit expression for the Lorentzian Cayley form is
\be\label{Phi-L-intr}
\Phi_L = \frac{\im}{2} \Omega_+ - \frac{\im}{2} \Omega_- + \frac{1}{2} \omega_r\wedge \omega_r.
\ee
The Lorentzian Cayley form (\ref{Phi-L-intr}) admits an expression that is adapted to the calibrations that it defines, and it is this expression that motivates our interest in this object. 

Before we discuss calibrations for $\Phi_L$ we would like to remind the reader the story for the real split Cayley form. This is a real 4-form that is built via (\ref{Phi-spinor}) from a real unit spinor of ${\rm Spin}(4,4)$. The stabiliser of the real split Cayley form is ${\rm Spin}(4,3)$, which is also the stabiliser of a real unit spinor. There are two different types of calibrations that this form defines. Thus, it has calibrations $H$ which are Riemannian in the restriction of the metric $g_\Phi$ to $H$. In a basis adapted to such a calibration the split Cayley form takes the form (\ref{Phi-split-calibr}). But it also has calibrations $H$ that have split signature metric. In a basis adapted to such a calibration the Cayley form is written as (\ref{Phi-split-calibr-2}). The two different types of calibrations in the split case can be understood as corresponding to a copy of $\Hq$ in the space of split octonions $\Oc'$, or to a choice of a copy of split quaternions $\Hq'\subset \Oc'$. 

The story for the Lorentzian Cayley form is somewhat analogous. As for the real split Cayley form, there are two different types of calibrations that $\Phi_L$ has. They can now be distinguished by whether the restriction of $\Phi_L$ to $H$ is purely real or purely imaginary. Thus, we will say that a four-dimensional subspace $H\subset TM$ is {\bf Lorentzian-calibrated} by $\Phi_L$ if for any orthonormal basis $\xi_{1,2,3,4}$ for $H$ we have $\Phi_L(\xi_1,\xi_2,\xi_3,\xi_4)=\pm \im$. Similarly, we will say that $H\subset TM$ is {\bf split-calibrated} by $\Phi_L$ if for any orthonormal basis $\xi_{1,2,3,4}$ for $H$ we have $\Phi_L(\xi_1,\xi_2,\xi_3,\xi_4)=\pm 1$. 
The Lorentzian Cayley form $\Phi_L$ can then be written in a form adapted to one of its calibrations. The most interesting for us expression is the one adapted to one of its Lorentzian calibrations. We have
\be\label{Lorentz-Cayley}
\Phi_L=
-\frac{1}{6} \Sigma^i_L \Sigma_L^i - \frac{1}{6}{\Sigma'_L}^i{\Sigma'_L}^i  - \Sigma_L^i {\Sigma'_L}^i ,
\ee
where $\Sigma_L^{1,2,3}$ and ${\Sigma'_L}^{1,2,3}$ are the Lorentzian (and thus complex) self-dual 2-forms for $H,H^\perp$. In particular, this formula makes it clear that the Lorentzian Cayley form is calibrated by copies of space with Lorentzian signature metric, and that in a basis adapted to such a calibration it is built from the complex self-dual 2-forms for $H,H^\perp$. This observation gives the justification for the title of this paper. While the stabiliser of $\Phi_L$ is ${\rm Spin}(3,3)={\rm SL}(4,\R)$, only the Lorentz ${\rm Spin}(3,1)$ subgroup of this is manifest in the way (\ref{Lorentz-Cayley}) of writing the complex 4-form $\Phi_L$. This is the diagonal Lorentz group that acts simultaneously by pseudo-orthogonal transformations on $H,H^\perp$.

Acting on the 4-form $\Phi_L$ with ${\rm GL}(8,\R)$ produces an orbit in the space of complex 4-forms of the algebraic type of $\Phi_L$. Any 4-form of this algebraic type defines a real metric of signature $(4,4)$ plus a complex unit semi-spinor $\psi_L$ of ${\rm Spin}(4,4)$ that has the property $\langle \hat{\psi}_L, \psi_L\rangle =0$. This can be confirmed by dimensions count. Indeed, we have
\be
\{ \text{Lorentzian Cayley forms}\} = {\rm GL}(8,\R)/{\rm Spin}(3,3).
\ee
The dimension of this space is $64-15=49$. This is the dimension $36$ of the space of metrics on $\R^8$ plus the dimension $16-2-1=13$ of the space of unit complex spinors that are orthogonal to their complex conjugates. So, we can write
\be\label{fibre-complex-intr}
\left\{ \lower1.2ex\vbox{\hbox{\text{Lorentzian}} \hbox{\text{Cayley forms}}}\right\}  = \{ \text{metrics on} \,\, \R^{4,4}\} \oplus \{ \text{complex unit spinors} \,\, \psi_L: \langle \psi_L,\hat{\psi}_L\rangle=0\}.
\ee

The definition of $\Phi_L$ and the expression (\ref{Lorentz-Cayley}) for the Lorentzian Cayley form solves a puzzle that existed in the context of geometry in four dimensions. Thus, it is well-known that a conformal metric on a four-dimensional manifold can be encoded into the knowledge of the Hodge operator on 2-forms, and thus in the knowledge of which 2-forms are self-dual. This works for all three different possible signatures in this number of dimensions, and can be described very concretely by means of the so-called Urbantke formula \cite{Urbantke:1984eb}. Let $B^i, i=1,2,3\in\Lambda^2(M)$ be a triple of 2-forms that is non-degenerate in the sense that the symmetric $3\times 3$ matrix of wedge-products $B^i\wedge B^j, i,j=1,2,3$ is non-singular. The triple then defines the metric 
\be\label{Urbantke}
g_B(\xi,\eta)v_B = \frac{1}{6} \epsilon^{ijk} \xi \lrcorner B^i \wedge \eta \lrcorner B^j \wedge B^k.
\ee
Here $\epsilon^{ijk}$ is the completely anti-symmetric rank 3 tensor. The right-hand side of this expression is a top degree form, and is a symmetric pairing of two vector fields $\xi,\eta\in TM$. It is required to be equal to a symmetric pairing $g_B(\xi,\eta)$ times the volume form $v_B$ of $g_B$. Thus, the formula (\ref{Urbantke}) defines a metric $g_B$ and, in particular, a conformal metric on $M$. The 2-forms $B^i$ are then self-dual (or anti-self dual, this depends on the orientation chosen) in the conformal metric $g_B$.  

The formula (\ref{Urbantke}) is a version of a statement well-known in the mathematical literature. For example, when we in addition demand that $B^i\wedge B^j\sim \delta^{ij}$ and that $dB^i=0$, we get the hyper-K\"ahler triples of Donaldson, see \cite{Donaldson}.

The formula (\ref{Urbantke}) works slightly differently for different signatures. When the 2-forms $B^i$ are real 2-forms on $M$, one can obtain via (\ref{Urbantke}) only the metrics of the Euclidean and split signatures. This is not surprising, because the 2-forms $B^i$ become identified with self-dual 2-forms of the metric $g_B$, and self-dual 2-forms are real only in the Euclidean and split signatures. These two signatures are distinguished by the properties of the matrix $B^i\wedge B^j$ of the wedge products of the 2-forms. This is a symmetric $3\times 3$ matrix that can be diagonalised by a ${\rm GL}(3)$ transformation. There are just two possibilities for the relative signs of the eigenvalues: either they are all of the same sign and the matrix arising is definite, or this matrix is indefinite. It is the first of these two cases for which the metric arising in (\ref{Urbantke}) is Euclidean. The indefinite case gives the split signature metric via (\ref{Urbantke}).

The Lorentzian case is more subtle. One must start with complex-valued 2-forms $B^i\in \Lambda^2(M,\C)$. In the Lorentzian signature the anti-self dual 2-forms are complex-conjugate of the self-dual ones. Since anti-self dual 2-forms are wedge product orthogonal to the self-dual ones, and we obtain the metric $g_B$ by declaring the 2-forms $B^i$ to be the self-dual 2-forms of $g_B$, we must impose the conditions
\be
B^i \wedge \overline{B^j}=0.
\ee
These are nine "reality conditions". They can be shown to guarantee that the conformal metric produced by (\ref{Urbantke}) from the complex 2-forms $B^i$ is real Lorentzian. The Lorentzian version of (\ref{Urbantke}) also requires an addition factor of the imaginary unit on the right-hand side. For more details on this story in different signatures we refer the reader to e.g. \cite{Krasnov:2020lku}, see Section 5.4. 

While the formula (\ref{Urbantke}) appears to be a curiosity of four-dimensional geometry, it is closely related to the geometry of Cayley forms in eight dimensions by a variant of dimensional reduction, at least in its Riemannian and split signature versions. Indeed, starting with a real Cayley form $\Phi$ of the Riemannian type, i.e. one that is compatible with a Riemannian signature metric, and choosing a calibrating subspace $H\subset TM$, we can choose a basis of self-dual bivectors $\sigma^{1,2,3}\in \Lambda^+ H$. Inserting these into $\Phi$ and restricting the resulting 2-forms onto $H^\perp$, we get a triple of 2-forms $B^i:= \sigma^i\lrcorner \Phi|_{H^\perp}$ on $H^\perp$. The conformal metric on $H^\perp$ defined by these $B^i$ via (\ref{Urbantke}) is then the same as the conformal metric induced on $H^\perp$ by $g_\Phi$. In this way (\ref{Urbantke}) can be understood as the dimensional reduction of the formula (\ref{metric-Phi}). 

This can also be phrased in terms of octonions. The Cayley form $\Phi$ of Riemannian type identifies the tangent space to every point with a copy of octonions $\Oc$. Choosing a copy of quaternions $\Hq\subset \Oc$ then gives the dimensional reduction to four dimensions, which is behind the Euclidean version of the formula (\ref{Urbantke}). Similarly, the real split Cayley 4-form $\Phi$ identifies $\R^8$ with the space of split octonions $\Oc'$. Choosing a copy $\Hq'\subset \Oc'$ gives the setup for dimensional reduction that is responsible for the split signature version of (\ref{Urbantke}). It is clear that there is no room here for a similar dimensional reduction interpretation of the Lorentzian version of the Urbantke formula, as long as one uses real Cayley forms. 

Nevertheless, the Lorentzian version of the Urbantke formula exists, and it is a natural question whether there is a version of the Cayley form that is responsible for the Lorentzian Urbantke formula, via dimensional reduction. The Lorentzian Cayley form (\ref{Lorentz-Cayley}) is exactly this object, and the main goal of this paper is to explain the geometry associated with such a form.

The organisation of this paper is as follows. We start, in Section \ref{sec:definitions}, by reviewing the linear algebra behind the notion of the Cayley form. We discuss the triple cross product and the notion of a compatible with a 4-form metric. A choice of a unit vector leads to cross-products and normed algebras. We also review how every Cayley form is given by a square of a certain unit Weyl spinor. We discuss Riemannian Cayley forms in more details in Section~\ref{sec:Riemannian}. We introduce an octonionic model of the Clifford algebra ${\rm Cl}(8)$, which then allows to do explicit calculations with Cayley forms corresponding to the same metric. We also review the notion of pure spinors here, and show how a complex viewpoint on a real Cayley form arises by representing a unit spinor as a sum of two complex conjugate pure spinors. We then introduce the notion of complex Cayley forms, and obtain an explicit expression for such a form in a basis adapted to a calibration. We prove that the ${\rm SU}(4)$ structure defined by a complex Cayley form is integrable if and only if this form is closed. We then proceed to analogous constructions in the split signature case in Section \ref{sec:split}. We treat the real split Cayley forms here, and explain the associated geometry. In particular, there are now two different possible viewpoints on the real split Cayley form, the complex and para-complex ones. Section \ref{sec:Lorentzian} is central to the paper. Here we discuss the two types of possible complex split Cayley forms, and define the notion of the Lorentzian Cayley form, as a special member of one of the families of complex split Cayley forms. We show how from the point of view of a $4+4$ split the Lorentzian Cayley form is built from the complex Lorentzian self-dual 2-forms. We also prove that the ${\rm SL}(4,\R)$ structure defined by a Lorentzian Cayley form is integrable if and only if this form is closed. The last section of the paper is an Appendix that explicitly computes the tangent space to the space of Lorentzian Cayley forms corresponding to the same metric.

\section{Cayley forms, cross products and normed algebras}
\label{sec:definitions}

We start with a series of definitions to remind the reader the notion of a Cayley form and its link to octonions. In this section we only consider the case of a positive definite metric. All definitions are from  \cite{Salamon:2010cs}. 

Let $W$ be a finite dimensional real vector space with inner product $W\times W\to \R$, which we denote by round brackets $(u,v)$. Thus, $(u,u):=|u|^2$ is the squared norm. 
\begin{definition} An alternating multi-linear map $W^3\to W (u,v,w)\to u\times v\times w$ is called a triple cross product if it satisfies $(u\times v\times w, u)=(u\times v\times w, v)=(u\times v\times w, w)=0$ and 
 \be
|u\times v\times w|^2 = |u\wedge v\wedge w|^2 := {\rm det}\left( \begin{array}{ccc} |u|^2 & (u,v) & (u,w) \\ (v,u) & |v|^2 & (v,w) \\ (w,u) & (w,v) & |w|^2\end{array}\right).
\ee
\end{definition}
Triple cross product exists only in dimensions $1,2,4,8$. It vanishes in dimensions $1,2$ and is unique in a given orientation in dimension 4. Let us now assume that the dimension is 8. In this dimension, there arises the notion of a Cayley form. 
\begin{definition} A 4-form $\Phi\in \Lambda^4 W^*$ is called non-degenerate if for any triple $u,v,w$ of linearly independent vectors in $W$ there exists a vector $x\in W$ such that $\Phi(u,v,w,x)\not=0$. An inner product in $W$ is called {\bf compatible} with a 4-form $\Phi$ if the map $W^3\to W (u,v,w)\to u\times v\times w$ defined by 
\be
( x, u\times v \times w) :=\Phi(x,u,v,w)
\ee
is a triple cross product. A 4-form $\Phi\in \Lambda^4 W^*$ is called a {\bf Cayley} form if it admits a compatible inner product. 
\end{definition}
An intrinsic characterisation of Cayley forms is provided by Theorem 7.8 in \cite{Salamon:2010cs}. If $\Phi$ is a Cayley form, then its compatible inner product is uniquely determined by $\Phi$. In particular, there is a unique orientation of $W$ and the associated volume form so that we have
\be
|u\wedge v|^2 \, {\rm vol} := ( |u|^2|v|^2- (u,v)^2) {\rm vol} = \frac{1}{6} u \lrcorner v \lrcorner \Phi \wedge u \lrcorner v \lrcorner \Phi \wedge \Phi.
\ee
The inner product defined by a Cayley form $\Phi$ can be explicitly extracted by the procedure explained in Lemma 7.13 of \cite{Salamon:2010cs}, see formula (7.20) of this reference. 

Choosing a unit vector $e\in W$ determines a skew-symmetric map $V_e\times V_e\to V_e$, where $V_e:= e^\perp$. This map is a cross product $u\times_e v:= u\times e \times v$. Let us discuss the notion of the cross-product. 

Let $V$ be a finite dimensional real vector space with inner product $V\times V\to \R$ that we denote by round brackets $(u,v)$. Then $(u,u):=|u|^2$ is the squared norm. 
\begin{definition}
A skew-symmetric bilinear map $V\times V\to V: u,v\to u\times v$ is called a cross-product if it satisfies (i) $(u\times v, u)=(u\times v, v)=0$ and (ii) $|u\times v|^2 = |u|^2 |v|^2- (u,v)^2$. 
\end{definition}
A cross-product exists only in dimensions 0,1,3,7, see Theorem 2.5 of \cite{Salamon:2010cs}. It is zero in dimensions 0,1, and in dimension 3 it is unique in a given orientation of $V$. The situation is much more interesting in dimension 7. A cross product on $V$ endows it with an alternating 3-form $\phi\in \Lambda^3 V^*$ (called the {\bf associative calibration}) defined via
\be\label{phi-cross-prod}
\phi(u,v,w):= (u\times v, w).
\ee

One can also start with a 3-form $\phi\in \Lambda^3 V^*$ and ask whether $\phi$ determines an inner product. 
\begin{definition} Let $V$ be a real vector space. A 3-form $\phi\in \Lambda V^*$ is called non-degenerate if, for every pair of linearly independent vectors $u,v\in V$ there exists a vector $w\in V$ such that $\phi(u,v,w)\not=0$. An inner product on $V$ is called compatible with $\phi$ if the map $\times: V\times V\to V$ defined by (\ref{phi-cross-prod}) is a cross product.
\end{definition}
The case of dimension seven is then particularly interesting because in this dimension a generic non-degenerate $\phi$ has a compatible inner product and moreover determines its compatible inner product uniquely. This is the content of Theorem 3.2 in \cite{Salamon:2010cs} that states that $\phi$ is non-degenerate if and only if it admits a compatible inner product, and moreover, this inner product if it exists is uniquely determined by $\phi$. Explicitly, given a non-degenerate $\phi \in \Lambda^3 V^*$, there exists an orientation of $V$ and the associated volume form ${\rm vol}\in \Lambda^7 V^*$ such that 
\be\label{phi-metric}
(u,v) {\rm vol} = \frac{1}{6} u\lrcorner \phi \wedge v \lrcorner \phi \wedge \phi .
\ee
A non-degenerate $\phi$ determines both the inner product in $V$ and the orientation ${\rm vol}$. Note, however, that for a non-degenerate real 3-form $\phi\in \Lambda^3 V^*$, the inner product resulting from (\ref{phi-metric}) may be both positive definite, as well as have signature $(4,3)$. There exists a degree seven invariant that can be constructed from $\phi$ which is non-zero if and only if $\phi$ is non-degenerate, and whose sign determines the signature of the inner product arising via (\ref{phi-metric}). This invariant is discussed, e.g. in \cite{Agricola}. 

Given a cross product one can consider the expression $(u\times v)\times w$. It can be seen that this is alternating on any triple of pairwise orthogonal vectors $u,v,w\in V$. It extends uniquely to an alternating map $V^3\to V$ called the associator bracket given by
\be
[u,v,w]:= \frac{1}{3}( (u\times v)\times w+(v\times w)\times u+(w\times u)\times v).
\ee
This gives rise to a 4-form $\psi\in \Lambda^4 V^*$ called the {\bf coassociative calibration} defined via
\be
\psi(u,v,w,x):= ( [u,v,w],x).
\ee
It can be shown that $\psi={}^* \phi$, where ${}^*$ is the Hodge operator ${}^*: \Lambda^k V^*\to \Lambda^{7-k} V^*$ associated to the inner product, and the orientation is as explained in Lemma 4.8 of \cite{Salamon:2010cs}.

Let us now discuss the relation between cross products, triple cross products, and normed algebras. The first statement is that a cross product determines a normed algebra and every normed algebra determines a cross product. 
\begin{definition} A normed algebra $W$ is a finite dimensional real vector space equipped with an inner product and a bilinear map $W\times W\to W$ (called the product), and also a unit $1\in W$. These data satisfy $1 u=u1=u$ and $|uv|=|u||v|$. 
\end{definition}
We have the notion of conjugation on $W$, defined by $\bar{1}=1$ and $\bar{u}=-u, u\in 1^\perp$. Any element of $W$ that is orthogonal to the unit element is called imaginary, and the conjugation changes the sign of the imaginary elements. The normed algebras and vector spaces equipped with cross products are related as follows, see Theorem 5.4 of \cite{Salamon:2010cs}. In one direction, if $W$ is a normed algebra then $V=1^\perp\subset W$ is equipped with the cross product defined as 
\be
u\times v = uv + (u,v), \qquad u,v\in 1^\perp.
\ee
In the opposite direction, if $V$ is a finite dimensional vector space equipped with an inner product and cross product, then $W=\R \oplus V$ is a normed algebra with the product
\be
uv = u_0 v_0 - (u_1, v_1) + u_0 v_1 + v_0 u_1 + u_1\times v_1,
\ee
where we write $u=u_0+u_1, v=v_0+v_1, u_0,v_0\in \R, u_1, v_1\in V$. The normed algebras one gets from the cross products in dimensions 0,1,3,7 are algebras isomorphic to $\R, \C, \Hq, \Oc$ respectively. 

Let us now explain how a Cayley form together with a choice of a unit vector $e\in W$ identifies $W$ with a copy of $\Oc$. This is the content of Theorem 6.15 of \cite{Salamon:2010cs}. 
\begin{theorem} Assume ${\rm dim}(W)=8$, let $\Phi$ be a Cayley form with its associated triple cross product, and let $e\in W$ be a unit vector. 
\begin{itemize}
\item Define the map $\psi_e: W^4\to \R$ by
\be
\psi_e(u,v,w,x) := ( e\times u\times v,e\times w\times x) - ( (u,w)-(u,e)(e,w))((v,x)-(v,e)(e,x)) \\ \nonumber
+ ( (u,x)-(u,e)(e,x))((v,w)-(v,e)(e,w)).
\ee
Then this map is alternating $\psi_e\in \Lambda^4 W^*$ and
\be
\Phi = e^* \wedge \phi_e + \epsilon \psi_e, \qquad \phi_e = e\lrcorner \Phi \in \Lambda^3 W^*,
\ee
and $\epsilon=\pm 1$ is a sign that is explained in Lemma 6.10 of \cite{Salamon:2010cs}.
\item The subspace $V_e:= e^\perp$ carries a cross product 
\be
V_e\times V_e\to V_e: u\times_e v := u\times e\times v.
\ee
The restriction of $\phi_e$ to $V_e$ is the associative calibration, and the restriction of $\psi_e$ to $V_e$ is the coassociative calibration. 
\item The space $W$ is a normed algebra with unit $e$ and multiplication and conjugations given by
\be
uv:= u\times e\times v+ (u,e)v + (v,e)u - (u,v)e, \qquad \bar{u} = 2(u,e) e - u.
\ee
\end{itemize}
\end{theorem}

The final statement that we need is Theorem 10.3 from \cite{Salamon:2010cs}. This shows how, given a Cayley form $\Phi\in \Lambda^4 W^*$, spinors for ${\rm Spin}(8)$ can be identified with certain subspaces in the spaces of differential forms. The Cayley form itself is given by the construction (\ref{Phi-spinor}) from a particular canonical spinor. More precisely, we have
\be
S^+ = \Lambda^0 \oplus \Lambda^2_7, \qquad S^- = \Lambda^1,
\ee
where
\be
\Lambda^2_7 := \{ \omega\in \Lambda^2 W^* : *(\Phi\wedge \omega) = 3\omega\} = \\ \nonumber
\{ u^*\wedge v^* - u\lrcorner v\lrcorner \Phi: u,v\in W\}.
\ee
The Clifford algebra is then given by
\be
\Gamma(u) = \left(\begin{array}{cc} 0 & \gamma(u)^* \\  \gamma(u) & 0 \end{array}\right),
\ee
where $\gamma(u): S^+\to S^-$ and is given explicitly by
\be
\gamma(u)(\lambda,w) = \lambda u^* + 2 u\lrcorner \omega, \qquad \lambda\in \Lambda^0, \omega\in \Lambda^2_7, \\ \nonumber
\gamma(u)^*(\theta) = ( u\lrcorner \theta, u^*\wedge \theta\Big|_{\Lambda^2_7}).
\ee
Then $\gamma(u)^* \gamma(u) = |u|^2 \id$, and we have the Clifford algebra $\Gamma(u)\Gamma(u)=|u|^2 \id$. Note that our Clifford algebra conventions are different from those in \cite{Salamon:2010cs}. The theorem 10.3 from \cite{Salamon:2010cs} then gives an explicit expression for the triple cross product that corresponds to the special spinor $\psi=(1,0)\in S^+$. The corresponding to it Cayley form is given by (\ref{Phi-spinor}) with spinor $\psi=(1,0)$. This shows that every Cayley form defines a metric and then an associated Clifford algebra, as well as a spinor, in which it then takes the form (\ref{Phi-spinor}). 

\section{Spinors of ${\rm Spin}(8)$ and Riemannian Cayley forms}
\label{sec:Riemannian}

For the purpose of doing explicit calculations with Cayley forms, it is very convenient to have an explicit model for the Clifford algebra ${\rm Cl}(8)$. In the previous section we explained how such a model arises when $\Phi$ is chosen. The Cayley form is then produced via (\ref{Phi-spinor}) from a particular canonical spinor that $\Phi$ defines. However, for practical calculations it is much more convenient to fix the metric and the corresponding Clifford algebra, together with a model for it, and only allow the spinor $\psi$ to vary. All calculations in the remainder of this paper are of this type. 

Given that a Cayley form $\Phi$ identifies $TM$ with a copy of $\Oc$, there results an octonionic model for ${\rm Cl}(8)$. This model is most convenient to explicit calculations. The purpose of this section is to describe the octonionic model for ${\rm Cl}(8)$, and then use it to derive some properties (and explicit expressions) for the Riemannian Cayley form $\Phi$. 

\subsection{Octonions}
 
 The space of octonions $\Oc$ is a normed algebra with the property $|xy|=|x| |y|$ (i.e. a composition algebra). The usual octonions (unlike split octonions) also have the property that the norm of every non-zero element is not zero, which makes them into a division algebra. It is non-commutative and non-associative, but alternative, which can be stated as the property that the subalgebra generated by any two elements is associative. 
 
 A general octonion is an object
 \be
 q= q_0 \id +\sum_{a=1}^7 q_a {\bf e}^a,
 \ee
 where ${\bf e}^a$ are unit imaginary octonions. The unit octonions anti-commute and square to minus the identity. The octonion conjugate changes the sign of all the imaginary octonions. The octonionic pairing is
 \be\label{oct-pairing}
 (q,q)=|q|^2= q\overline{q} = (q_0)^2 + \sum_{a=1}^7 (q_a)^2.
 \ee
 
 We encode the octonionic product by the cross-product in the space of imaginary octonions. Thus, we write
 \be
 \Oc = \R + {\rm Im}\,\Oc.
 \ee
 Let $\e^{1,\ldots,7}$ be a basis in the space of imaginary octonions. The cross-product in ${\rm Im}\,\Oc$ is then encoded by the associator 3-form on $\R^7$, which we chose to be
 \be\label{3-form}
\phi = e^{567} + e^5( e^{41}-e^{23}) + e^6(e^{42}-e^{31}) + e^7(e^{43}-e^{12}).
 \ee
Here $e^a, a=1,\ldots,7$ is a basis of one-forms on $\R^7$, and the notation is $e^{abc\ldots}=e^a\wedge e^b\wedge e^c \ldots$. This encodes the cross-product in the sense that
\be
( \e^a\times \e^b, \e^c)=(e^c)^* \lrcorner (e^b)^* \lrcorner (e^a)^* \lrcorner  \phi, 
\ee
where $(e^a)^*$ is the dual basis of vectors. So, for instance $\e^5 \times \e^6= \e^7$. For later purposes we will also mention that the dual 4-form is given by
 \be\label{dual-C}
 {}^* \phi = e^{1234} +e^{67}( e^{41}-e^{23}) + e^{75}(e^{42}-e^{31}) + e^{56} (e^{43}-e^{12}).
 \ee

 \subsection{Octonionic model for ${\rm Cl}(8)$}
 
 The Clifford algebra ${\rm Cl}(8)$ can be generated by the operators of, say, left multiplication by unit octonions. Concretely, we have the following representation of a general linear combination of the Clifford generators
 \be\label{cl8-oct}
 q_0 \Gamma_0 + \sum_{a=1}^7 q_a \Gamma_a = \left( \begin{array}{cc} 0 & L_{\bar{q}} \\ L_q & 0 \end{array}\right). 
 \ee
Thus, ${\rm Cl}(8)$ is generated by $2\times 2$ matrices with octonionic entries that act on 2-component columns with entries in $\Oc$
\be
\psi = \left(\begin{array}{c} \alpha \\ \beta \end{array}\right), \qquad \alpha,\beta \in \Oc.
\ee
Majorana-Weyl spinors are then identified with copies of octonions $S^\pm= \Oc$. We will take the Majorana-Weyl spinors in $S^+$ to be those with the upper element different from zero, i.e. $\alpha\not=0, \beta=0$. The ${\rm Spin}(8)$ invariant inner product on $S^+$ is just the pairing (\ref{oct-pairing}) of two octonions. 

For concrete computations, it is convenient to have an explicit form of the $\Gamma$-matrices. To this end we arrange octonions into 8-component columns that start with $q_0$ as the top entry. In this representation the $\Gamma$-matrices (\ref{cl8-oct}) become the following $16\times 16$ matrices
\begin{equation}\label{Gamma-matr-spin8}
    \Gamma_0=
    \begin{pmatrix}
    0&\id \\
    \id &0
    \end{pmatrix}
    , \quad
    \Gamma_a=
    \begin{pmatrix}
    0&-E_a\\
    E_a&0
    \end{pmatrix}
    \ \text{for} \ a\in \{1,\ldots,7\},
\end{equation}
where
\begin{equation}\label{E-matr-spin8}
    \begin{split} 
        E_1&=-E_{01}+E_{27}-E_{36}+E_{45},\\
        E_2&=-E_{02}-E_{17}+E_{35}+E_{46},\\
        E_3&=-E_{03}+E_{16}-E_{25}+E_{47},\\
        E_4&=-E_{04}-E_{15}-E_{26}-E_{37},\\
        E_5&=-E_{05}+E_{14}+E_{23}-E_{67},\\
        E_6&=-E_{06}-E_{13}+E_{24}+E_{57},\\
        E_7&=-E_{07}+E_{12}+E_{34}-E_{56},
    \end{split}
\end{equation}
and $E_{ij}$ are the usual generators of $\mathfrak{so}(8)$, i.e. anti-symmetric matrices with $+1$ on the $i$th row and $j$the column. 

\subsection{Cayley form}

We now take a Majorana-Weyl spinor $\psi$ given by $\alpha=\id,\beta=0$. This is a unit spinor $\langle \id,\id\rangle=1$. The insertion of two distinct $\Gamma$-matrices between this spinor vanishes $\langle \id, \Gamma\Gamma \id\rangle=0$. The insertion of four distinct $\Gamma$-matrices gives the following 4-form on $\R^8$
\be\label{Cayley-plus}
\Phi:=\langle \id, \Gamma\Gamma\Gamma\Gamma \id \rangle = e^0 \wedge \phi -{}^* \phi, 
\ee
where $\phi$ is given by (\ref{3-form}) and ${}^*\phi$ is its dual 4-form given by (\ref{dual-C}). The meaning of the above equation is that it lists non-zero values for the insertion of four different $\Gamma$-matrices between two copies of the spinor $\alpha=\id,\beta=0$. Thus, e.g. $(\id, E_1 E_2 E_3 E_4 \id)=-1$.
 
\subsection{Cayley form for an arbitrary unit spinor}

For later purposes, we give the formulas for the Cayley form obtained as 
\be
\Phi_\psi:=\langle \psi, \Gamma\Gamma\Gamma\Gamma \psi \rangle,
\ee
where $\psi$ is the unit spinor $\langle\psi,\psi\rangle=1$ parametrised as
\be
\psi = ( \sqrt{1-|\alpha|^2}, \alpha), \qquad \alpha\in \R^7.
\ee
We have
\be\label{fibre-unit-spinor}
\Phi_\psi = e^0 \wedge C_\psi - {}^* C_\psi,
\ee
where
\be
C_\psi := (1-2|\alpha|^2) C + 2\alpha \wedge \alpha^* \lrcorner C - 2\sqrt{1-|\alpha|^2} \alpha^* \lrcorner {}^*C, \\ 
{}^* C_\psi = {}^*C - 2\alpha \wedge \alpha^* \lrcorner {}^* C + 2\sqrt{1-|\alpha|^2} \alpha\wedge C,
\ee
where $\alpha$ is interpreted as a 1-form in $\R^7$ in these formulas, and $\alpha^*$ is the dual vector field, computed with respect to the standard metric on $\R^7$. 

The formula (\ref{fibre-unit-spinor}) gives an explicit description of the 4-forms of the Cayley algebraic type that give rise to the same metric. It can also be noted that the tangent space to the space of such 4-forms, at $\Phi$, is given by the 4-forms of the type
\be\label{delta-Phi}
\delta \Phi = e^0 \wedge Y^* \lrcorner \Phi-  Y\wedge (e^0)^* \lrcorner \Phi,
\ee
where $Y\in \R^7$. 
 
\subsection{Pure spinors}

We will now remind basic facts about pure spinors. Let $V$ be a real vector space of dimension ${\rm dim}(V)=2n$, equipped with a metric that is not necessarily definite. Given a spinor $\psi$ we define $M(\psi):= V^0_\psi$, the subspace of the complexification $V_\C$ of $V$ that annihilates $\psi$ via the Clifford product. The spinor $\psi$ is said to be pure if the dimension of the space $M(\psi)$ is maximal possible, i.e. $n$. We shall refer to $M(\psi)$ as the {\bf maximal totally null} (MTN) subspace corresponding to the pure spinor $\psi$. It is known since Cartan \cite{Cartan} that pure spinors are Weyl. Cartan also gives a very useful algebraic characterisation of pure spinors, which we now explain.

Given a (Weyl) spinor $\psi$ (or a pair of Weyl spinors $\psi,\phi$), one can insert a number of $\Gamma$-matrices between two copies of $\psi$ (or, more generally, between $\psi$ and $\phi$). Let us introduce the convenient notation
\be
\Lambda^k_\C(\R^{2n})\ni B_k(\psi,\phi) := \langle \psi, \underbrace{\Gamma\ldots \Gamma}_{\text{$k$ times}} \phi\rangle.
\ee
It is assumed that distinct $\Gamma$-matrices are inserted, thus giving components of a degree $k$ differential form (in general complex) in $\R^{2n}$. 

The following proposition is due to Cartan \cite{Cartan} and Chevalley \cite{Chevalley}
\begin{proposition} \label{thm:cartan} A Weyl spinor $\phi$ is a pure (or in Cartan's terminology simple) if and only if $B_k(\phi,\phi) =0$ for $k\not=n$. The $n$-vector $B_n(\phi,\phi)$ is then proportional to the wedge product of the vectors constituting a basis of $M(\phi)$, where $M(\phi)$ is the MTN subspace that corresponds to $\phi$. 
\end{proposition}
We note that this theorem gives a practical way of recovering $M(\phi)$ as the set of vectors whose insertion into $B_n(\phi,\phi)$ vanishes. Thus, this theorem establishes a one-to-one correspondence between pure spinors $\phi$ (up to rescaling) and MTN subspaces $M(\phi)$. This theorem also gives a set of quadratic constraints that each pure spinor must satisfy. 

For a general signature there are several distinct types of pure spinors of ${\rm Spin}(k,l)$, see \cite{KT}. For our purposes in this section we only need the case of ${\rm Spin}(8)$. In this case there is only one type of pure spinors and it can be shown that ${\rm Spin}(2n)$ acts transitively on the orbits of pure spinors of fixed length $\langle\hat{\psi},\psi\rangle$ in $S^\pm$, where $\hat{\psi}$ is an appropriate complex conjugate of $\psi$. See e.g. \cite{Harvey} for a proof of the transitivity of the group action on the pure spinor orbit. 

We will also need statements about the possible ways that the MTN subspaces of pure spinors can intersect. Let $\phi,\psi$ being two pure spinors of ${\rm Cl}(2n)$ of the same parity, with the corresponding MTN subspaces $M(\phi), M(\psi)$. We have the following set of statements, most due to Cartan \cite{Cartan}
\begin{proposition}\label{thm:same} The MTN subspaces of two pure spinors of the same parity $\phi,\psi\in S_+$ can intersect along $n-2m, m\in\mathbb{N}$ null directions. The intersection along $n-2m$ directions occurs if and only if $B_k(\phi,\psi)=0$ for $k<n-2m$ and $B_{n-2m}(\phi,\psi)$ is non-zero. It is then decomposable and  its factors are the common null directions in $M(\phi), M(\psi)$. Two MTN's $M(\phi), M(\psi)$ do not intersect if and only if $B_0(\phi,\psi)\not=0$. 
\end{proposition}
For a proof, see e.g. \cite{Budinich:1989bg}, Proposition 9. We also have the following fact due to Chevalley
\begin{proposition} When the MTN subspaces of two pure spinors $\phi,\psi$ intersect along $n-2$ directions, the sum $\phi+\psi$ is a pure spinor. 
\end{proposition}
This last fact follows from the fact that all Weyl spinors in four dimensions are pure. Its proof can also be found in \cite{Budinich:1989bg}.

\subsection{Pure spinors of ${\rm Spin}(2n)$ and complex structures}

There is a bijective relation between the (projective, i.e. up to complex rescalings) pure spinors of ${\rm Spin}(2n)$ and orthogonal complex structures on $\R^{2n}$ consistent with a given orientation. We review it only briefly, for more information see e.g. \cite{Salamon} or \cite{Lawson:1998yr} Chapter IV, Section 9. 

In one direction, as we have seen, the pure spinor $\psi\in S_+$ comes with its MTN subspace $M(\psi)$. There is a ${\rm Spin}(2n)$ invariant operation of complex conjugation $\psi\to \hat{\psi}$, which is anti-linear and squares to plus or minus the identity, depending on $n$. Depending on $n$, the complex conjugation operator either preserves $S_+$, or maps $S_+\to S_-$, see e.g. \cite{Bhoja:2022dsm} for an explicit description. The annihilator subspace $M(\hat{\psi})$ is the complex conjugate of $M(\psi)$, and $M(\psi),M(\hat{\psi})$ are complementary in $\R^{2n}_\C$ and together span it. One can define an orthogonal complex structure (OCS) on $\R^{2n}$ by declaring $M(\psi)$ to be the eigenspace of $J$ of eigenvalue $+\im$, and then $M(\hat{\psi})$ to be the eigenspace of eigenvalue $-\im$. This defines an OCS $J$, which in turn gives $\R^{2n}$ an orientation. 

To establish the relation in the other direction, one starts with an OCS $J$ on $\R^{2n}$. This defines the eigenspaces $E^\pm$
\be
J E^\pm = \pm \im E^\pm
\ee
of $J$ that are totally null, complex conjugate of each other and of dimension $n$. One then realises ${\rm Cl}(2n)$ in terms of creation-annihilation operators on $\Lambda(E^+)$, as is reviewed in e.g. \cite{Bhoja:2022dsm}. Then the top and  the zero degree elements of $\Lambda(E^+)$ are pure spinors that correspond to $J$. They may both be in $S_+$, or one in $S_+$ and one in $S_-$, depending on the dimension. To obtain pure spinors that correspond to other OSC one uses the action of ${\rm Spin}(2n)$. This action is transitive on both the space of OCS on $\R^{2n}$ that are compatible with a given orientation and on the space of pure spinors (of a given parity). This establishes that a given OSC can be linked to its corresponding pure spinor. 

We note that there are always two distinct pure spinors (up to complex multiplication) that correspond to a given OCS $J$. The MTN of one of these pure spinors is $E^+$, while the MTN of the other is $E^-$.

\subsection{Pure and non-pure spinors of ${\rm Spin}(8)$}

The octonionic model of ${\rm Cl}(8)$ identifies Weyl spinors of ${\rm Spin}(8)$ with complexified octonions. The only spinor bilinears $B_k(\psi,\psi)$ that exist are $B_0, B_4$. Cartan's purity criterion tells us that for a pure spinor $B_0(\phi,\phi)=0$, and this is the only condition that the pure spinor must satisfy. Thus, pure spinors of ${\rm Spin}(8)$ are complexified null octonions. Given the reviewed above relation between OCS compatible with a given orientation and projective pure spinors of a given parity, we see that the space $Z^+_8={\rm SO}(8)/{\rm U}(4)$ of OSC on $\R^8$ compatible with a given orientation can be realised as the null quadric $Q^6\subset \mathbb{P}^7$, see \cite{Salamon} for more on this interpretation. 

It is then easy to understand the other possible types of spinors that can exist in this dimension. Indeed,  any non-null complexified octonion can be written as a sum of two null complexified octonions. To see this, it is convenient to first rescale the complexified octonion by a complex number to make it have norm one. After this is done, we have $\alpha_1+\im\alpha_2$ with $|\alpha_1|^2-|\alpha_2|^2=1$ and $(\alpha_1,\alpha_2)=0$. We can then define
\be
|\alpha_1|=\cosh\tau, \qquad |\alpha_2|= \sinh\tau,
\ee
so that 
\be\label{compl-psi}
\alpha_1+\im \alpha_2 =  \cosh\tau \frac{\alpha_1}{|\alpha_1|} + \im \sinh\tau \frac{\alpha_2}{|\alpha_2|}
= \frac{1}{2}e^\tau \left(\frac{\alpha_1}{|\alpha_1|} + \im \frac{\alpha_2}{|\alpha_2|}\right) +
\frac{1}{2}e^{-\tau}  \left(\frac{\alpha_1}{|\alpha_1|} - \im \frac{\alpha_2}{|\alpha_2|}\right).
\ee
This means that any unit impure spinor of ${\rm Cl}(8)$ can be written as a linear combination of two complex conjugate pure spinors. The same applies for the case of a real octonion $\alpha_1$. In this case we just need to select an (arbitrary) octonion $\alpha_2$ orthogonal to $\alpha_1$ and of the same norm. Then 
\be
\alpha_1= \frac{1}{2}(\alpha_1+\im \alpha_2) + \frac{1}{2}(\alpha_1-\im\alpha_2)
\ee
is a sum of two null complexified octonions. 

Given that we need at most two pure spinors to represent an arbitrary spinor in eight dimensions, we can understand the geometry of impure spinors. The MTN of a pure spinor in eight dimensions is four-dimensional. If we take two pure spinors $\phi,\psi$ of the same parity, their MTN's can intersect in four, two and zero dimensions. When they intersect in four dimensions $\phi$ is proportional to $\psi$. When they intersect in two dimensions, the sum is still a pure spinor. When they intersect in zero dimensions, the sum is impure. Thus, any impure spinor in eight dimensions is a sum of two pure spinors whose MTN's do not intersect. In terms of complexified octonions, the impure spinors are given by a linear combination of two null octonions that are complex conjugate to each other. 

The case of complex impure spinors is then different from the real case in terms of the geometry that arises. In the complex case the stabiliser of $\phi=\alpha_1+\im \alpha_2$ must stabilise both $\alpha_{1,2}$. It is then a copy of ${\rm Spin}(6)={\rm SU}(4)$ that stabilises the two complex conjugate pure spinors composing $\phi$. When the spinor is real $\phi=\alpha_1$, the stabiliser is a copy of ${\rm Spin}(7)$ that stabilises the octonion $\alpha_1$. 

\subsection{Complex basis for the Cayley form}

Cayley form on $\R^8$ is a 4-form of a special algebraic type that endows $\R^8$ with a metric and a unit real Weyl spinor. As discussed in a previous subsection, we can always write a real unit spinor as the sum of a pure spinor and its complex conjugate. There is of course no unique way of doing this, and the choice involved is the same as the choice of a complex structure on $\R^8$. 

Thus, let us choose an orthogonal complex structure on $\R^8$. There is the corresponding pair of complex conjugate pure spinors $\psi_p,\hat{\psi}_p$, and their sum
\be\label{psi-real}
\psi = \psi_p +\hat{\psi}_p
\ee
is a real spinor. We rescale $\psi$ to make it unit, which is equivalent to $\langle \hat{\psi}_p,\psi \rangle =1/2$. The Cayley form for $\psi$ can then be computed as 
\be\label{cayley-comput}
\langle \psi, \Gamma\Gamma\Gamma\Gamma \psi \rangle =\langle \psi_p, \Gamma\Gamma\Gamma\Gamma \psi_p \rangle + \langle \hat{\psi}_p, \Gamma\Gamma\Gamma\Gamma \hat{\psi}_p \rangle + 2\langle \hat{\psi}_p, \Gamma\Gamma\Gamma\Gamma \psi_p \rangle .
\ee
The first term here is decomposable and is given by the product of the null directions of $\psi_p$. In other words, it is a multiple of the $(4,0)$ holomorphic top form $\Omega$ for the complex structure $J$
\be
\langle \psi_p, \Gamma\Gamma\Gamma\Gamma \psi_p \rangle = \frac{1}{2} \Omega.
\ee
Thus, the first two terms in (\ref{cayley-comput}) can be written as ${\rm Re}(\Omega)$. To understand the last term we note that
\be
\langle \hat{\psi}_p, \Gamma\Gamma \psi_p \rangle = \frac{\im}{2} \omega,
\ee
where $\omega$ is the $(1,1)$ Kahler form for $J$. The last term in (\ref{cayley-comput}) is then a multiple of $\omega\wedge \omega$
\be
\langle \hat{\psi}_p, \Gamma\Gamma \Gamma\Gamma \psi_p \rangle = -\frac{1}{4} \omega\wedge \omega.
\ee
Overall, we get
\be\label{Phi-compl-form}
\langle \psi, \Gamma\Gamma\Gamma\Gamma \psi \rangle = {\rm Re}(\Omega)-\frac{1}{2} \omega\wedge \omega.
\ee

To make this concrete, we take the spinor $\psi$ to be the identity octonion $\psi=\id$, and the pure spinor to be
\be
\psi_p = \frac{1}{2}(\id + \im \e^4).
\ee
If $\psi_p$ is interpreted as a positive parity spinor $\psi_p\in S_+$, and $\R^8$ is identified with octonions, the corresponding complex structure $J$ is $J=R_{\e^4}$, the operator of right multiplication by $\e^4$. The corresponding $+\im$ eigenspace is
\be
E^+ ={\rm Span}( e^4+ \im e^0, e^1+\im e^5, e^2+\im e^6, e^3+\im e^7). 
\ee
Then 
\be\label{Omega-omega-standard}
\Omega = (e^4+\im e^0)(e^1+\im e^5)(e^2+\im e^6)(e^3+\im e^7), \\ \nonumber
\omega = e^{15}+e^{26}+e^{37}+e^{40},
\ee
and
\be\label{Phi-CS}
\Phi = {\rm Re}(\Omega)-\frac{1}{2} \omega\wedge \omega= e^0 \wedge C -{}^*C,
\ee
where $C,{}^*C$ are given by (\ref{3-form}) and (\ref{dual-C}) respectively. 

Thus, every Cayley form is of the form $\langle \psi, \Gamma\Gamma\Gamma\Gamma \psi \rangle$, and can be written in the complex basis as (\ref{Phi-compl-form}) by choosing (non-uniquely) a complex structure on $\R^8$ whose associated pure spinor has $\psi$ as its real part. The stabiliser of $\Phi$ is the stabiliser of the spinor $\psi$, which is ${\rm Spin}(7)$. The stabiliser of a pure spinor $\psi_p$ whose real part is $\psi$ is a copy of ${\rm SU}(4)={\rm Spin}(6)$. Thus, the complex form (\ref{Phi-compl-form}) makes manifest only the ${\rm SU}(4)$ subgroup of its stabiliser subgroup ${\rm Spin}(7)$.

\subsection{Calibrations}

The expression (\ref{Phi-CS}) for the Cayley form $\Phi$ exhibits the ${\rm SU}(4)$ subgroup of its stabiliser group ${\rm Spin}(7)$. The same 4-form can be written differently, exhibiting a ${\rm SU}(2)\times{\rm SU}(2)\times{\rm SU}(2)$ subgroup. This also leads to the notion of 4-dimensional submanifolds of $\R^8$ calibrated by $\Phi$. 

We start by introducing a basis in the space of self-dual two-forms in the $e^{1,2,3,4}$ directions. Thus, we define
\be
\Sigma^1 = e^{41}- e^{23}, \quad \Sigma^2= e^{42}-e^{31}, \quad \Sigma^3= e^{43}-e^{12}.
\ee
We then define a similar basis for the directions $e^{0,5,6,7}$
\be
{\Sigma'}^1= e^{05}-e^{67}, \quad {\Sigma'}^1= e^{06}-e^{75}, \quad {\Sigma'}^1= e^{07}-e^{56}.
\ee
We note that the volume forms for each of these 4-dimensional subspaces can be written as
\be
e^{4123}= - \frac{1}{6} \Sigma^i \Sigma^i, \qquad e^{0567}= - \frac{1}{6} {\Sigma'}^i {\Sigma'}^i.
\ee
It is then not hard to check that the Cayley form (\ref{Cayley-plus}), which also equals (\ref{Phi-CS}), 
can be written as
\be\label{Phi-calibr}
\Phi = - \frac{1}{6} \Sigma^i \Sigma^i- \frac{1}{6} {\Sigma'}^i {\Sigma'}^i+ \Sigma^i{\Sigma'}^i.
\ee
The ${\rm SU}(2)\times{\rm SU}(2)\times{\rm SU}(2)$ invariance of this form is manifest. Indeed, one can perform the ${\rm SO}(4)$ rotations in each of the subspaces $e^{1,2,3,4}$ and $e^{0,5,6,7}$. The orthogonal group in four dimensions splits as ${\rm SO}(4)= {\rm SU}(2)\times{\rm SU}(2)/\Z_2$. One of these factors acts on the self-dual 2-forms $\Sigma^i$ as an ${\rm SO}(3)$ transformation, the other leaves it invariant. So, the subgroups ${\rm SU}(2), {\rm SU}'(2)$ that leave invariant the self-dual 2-forms $\Sigma^i, {\Sigma}'^i$ are manifest in (\ref{Phi-calibr}). Yet another ${\rm SU}(2)$ is the diagonal in the ${\rm SU}(2)\times{\rm SU}(2)$ subgroup where each factor acts on the self-dual 2-forms. Only the diagonal ${\rm SU}(2)$ survives as the invariance of the last term in (\ref{Phi-calibr}).

Cayley 4-form in $\R^8$ provides an example of calibrations, see \cite{Harvey:1982xk} and \cite{Harvey}. A 4-dimensional subspace $X$ of $\R^8$ is called a calibration if the restriction of $\Phi$ to $X$ coincides with the volume form of $X$ at each point. Thus, the local model for a calibrated subspace $X$ is that of the $e^{0,5,6,7}$ subspace of $\R^8$. 

We end with the statement describing the Grassmanian of Cayley 4-planes in $\R^8$ as a homogeneous group space. 
\begin{proposition} The group ${\rm Spin}(7)$ acts transitively on the Grassmanian $G(\Phi)$ of Cayley 4-planes in $\R^8$. The isotropy group is $K= {\rm SU}(2)\times{\rm SU}(2)\times{\rm SU}(2)/\Z_2$. Thus, $G(\Phi)={\rm Spin}(7)/K$. 
\end{proposition}
This is theorem 1.38 from \cite{Harvey:1982xk}.

\subsection{Complex Cayley form}

Having understood the real Riemannian Cayley form, as one produced as $\langle \psi, \Gamma\Gamma\Gamma\Gamma \psi\rangle$ from a real unit spinor $\psi$, we can describe its complex generalisation. This is produced in the same way, but now from a {\bf complex} unit spinor. It is clear that the metric defined by such a complex Cayley form will still be real. Indeed, the only fact that is used in the algebra of obtaining the metric $g_\Phi$ from the real Cayley form $\Phi= \langle \psi, \Gamma\Gamma\Gamma\Gamma \psi\rangle$ is the unity of the spinor $\psi$. And so the metric defined by $\Phi$ continues to be real even when $\psi$ is a complex unit spinor. For the complex Cayley form given by (\ref{Phi-c}) this can be checked by an explicit computation. 

We have already seen in (\ref{compl-psi}) that an arbitrary complex unit spinor can be written as a linear combination of two complex conjugate pure spinors 
\be \label{psi-tau}
\psi_\tau:= e^\tau \psi_p + e^{-\tau} \hat{\psi}_p,
\ee
where
\be
\psi_p := \frac{1}{2} \left( \frac{\alpha_1}{|\alpha_1|} + \im \frac{\alpha_2}{|\alpha_2|} \right),
\ee
and $\hat{\psi}_p$ is the complex conjugate complexified octonion, which is also null. When $\tau=0$ we get the real unit spinor that we considered in (\ref{psi-real}). 

It is now not difficult to compute the complex Cayley form 
\be
\Phi_\tau:= \langle \psi_\tau, \Gamma\Gamma\Gamma\Gamma \psi_\tau\rangle.
\ee
It is given by the expression (\ref{Phi-c}) quoted in the Introduction. Unlike $\Phi=\Phi_0$, the complex Cayley form is only invariant under ${\rm SU}(4)$, which is the stabiliser of $\psi_\tau$. The expression (\ref{Phi-c}) makes manifest the full stabiliser subgroup. 

\subsection{Complex Cayley form in a basis adapted to a calibration}

The appropriate notion of a calibration for a complex Cayley form is as follows. We say that $H\subset \R^8$ is calibrated by $\Phi_\tau$ if for any orthonormal basis $\xi_{1,2,3,4}\in H$ we have $\Phi_\tau(\xi_1, \xi_2,\xi_3,\xi_4) = \pm \cosh(2\tau)$. 

Having chosen a calibration, we can rewrite $\Phi_\tau$ in a basis adapted to its calibration. We do this for our standard form (\ref{Phi-c})  with $\Omega,\omega$ given by (\ref{Omega-omega-standard}). To this end it is most convenient to first rewrite $\Omega,\omega$ in a calibration adapted basis. And so we order the basis of orthonormal vectors in $\R^8$ so that it is spanned by $e^I, e'^I, I=1,2,3,4$, where
\be
e'^1=e^5, \quad e'^2=e^6, \quad e'^3=e^7, \quad e'^4=e^0.
\ee
We choose the orientation to be $\epsilon^{4123}=+1$. We can then write
\be\label{Omega-eep}
{\rm Re}(\Omega)= \frac{1}{24} \epsilon_{IJKL} e^{I} e^{J} e^{K} e^{L}+ \frac{1}{24} \epsilon_{IJKL} e'^{I} e'^{J} e'^{K} e'^{L} - \frac{1}{4} \epsilon_{IJKL} e^{I} e^{J} e'^{K} e'^{L}, \\ \nonumber
{\rm Im}(\Omega)= \frac{1}{6} \epsilon_{IJKL} e'^{I} e^{J} e^{K} e^{L} - \frac{1}{6} \epsilon_{IJKL} e^I e'^{J} e'^{K} e'^{L} , \\ \nonumber
\frac{1}{2}\omega\wedge \omega = -\frac{1}{4} e^I e^J e'^K e'^L \delta_{IK} \delta_{JL}.
\ee
 This form of writing these 4-forms makes manifest the diagonal ${\rm SO}(4)\subset{\rm Spin}(8)$ subgroup of the stabiliser ${\rm SU}(4)\subset{\rm Spin}(8)$.

In terms of this basis the real Cayley form is given by
\be
\Phi = {\rm Re}(\Omega) - \frac{1}{2}\omega\wedge \omega = \frac{1}{24} \epsilon_{IJKL} e^{I} e^{J} e^{K} e^{L}+ \frac{1}{24} \epsilon_{IJKL} e'^{I} e'^{J} e'^{K} e'^{L} \\ \nonumber
+ \frac{1}{4}( \delta_{I[K} \delta_{L]J} - \epsilon_{IJKL}) e^{I} e^{J} e'^{K} e'^{L}.
\ee
This can be written in terms of the self-dual 2-forms 
\be\label{Sigma-SD}
\Sigma^{IJ} := \frac{1}{2}( \delta_{I[K} \delta_{L]J} - \epsilon_{IJKL}) e^{K} e^{L}, \\ \nonumber
\Sigma'^{IJ} := \frac{1}{2}( \delta_{I[K} \delta_{L]J} - \epsilon_{IJKL}) e'^{K} e'^{L},
\ee
and takes the form (\ref{Phi-calibr}).

On the other hand, using the formulas (\ref{Omega-eep}) we get for $\Phi_\tau$
\be\label{Phi-tau-calibr}
 {\rm Re}(\Phi_\tau)= \frac{\cosh(2\tau)}{24} \epsilon_{IJKL} e^{I} e^{J} e^{K} e^{L}+ \frac{\cosh(2\tau)}{24} \epsilon_{IJKL} e'^{I} e'^{J} e'^{K} e'^{L} \\ \nonumber
 +  \frac{1}{4}( \delta_{I[K} \delta_{L]J} - \cosh(2\tau) \epsilon_{IJKL}) e^{I} e^{J} e'^{K} e'^{L}, \\ \nonumber
 {\rm Im}(\Phi_\tau)= \im \frac{\sinh(2\tau)}{6} \epsilon_{IJKL} e'^{I} e^{J} e^{K} e^{L} - \im \frac{\sinh(2\tau)}{6} \epsilon_{IJKL} e^I e'^{J} e'^{K} e'^{L} .
\ee
This can no longer be written in terms of self-dual 2-forms (\ref{Sigma-SD}) for $H,H^\perp$. This way of writing manifests only the ${\rm SO}(4)$ subgroup of the stabiliser ${\rm SU}(4)$ of the real and imaginary parts of $\Phi_\tau$.

\subsection{Complex Cayley forms and integrable ${\rm SU}(4)$ structures}

The introduced complex Cayley form is an interesting geometric object because it encodes the data $(g,J, \Omega)$ of a Riemannian metric, a compatible with it complex structure $J$ as well as a top holomorphic form $\Omega$ in a single geometric object $\Phi_\tau$. Thus, for any $\tau\not=0$ the complex Cayley form $\Phi_\tau$ can be called an ${\rm SU}(4)$ structure. It is easy to see that this ${\rm SU}(4)$ structure is integrable if and only of $d\Phi_\tau=0$. 

Indeed, in one direction, when the ${\rm SU}(4)$ structure is integrable we have $\nabla \omega=0, \nabla\Omega=0$, where $\nabla$ is the covariant derivative. This implies $d\omega=0, d\Omega=0$, and so also $d\Phi_\tau=0$. In the other direction when $d\Phi_\tau=0$ the imaginary part of this equation gives $d{\rm Im}(\Omega)=0$. However, with $\Omega\in \Lambda^{4,0}$, we have $d\Omega \in \Lambda^{4,1}\oplus \Lambda^{3,2}$, where we are using the notation of Chapter 3 of \cite{Salamon-book}. This means that projecting $d{\rm Im}(\Omega)=0$ on its irreducible components gives $d\Omega=0$. It is easy to see that this implies the integrability of the complex structure defined by $\Omega$. Indeed, this implies that for any $\theta\in \Lambda^{1,0}$ we have $d\theta|_{(0,2)}=0$, which is equivalent to the vanishing of the Nijenhuis tensor. 

The intrinsic torsion of an ${\rm SU}(4)$ structure, see \cite{Salamon-Chiossi} for a similar discussion in the case of ${\rm SU}(3)$ structure, takes values in the space
\be
\Lambda^1 \otimes {\mathfrak su}(4)^\perp = [[\Lambda^{1,0} \otimes \Lambda^{2,0}]] \oplus [[\Lambda^{2,1}]] \oplus \Lambda^1,
\ee
where the space $[[\Lambda^{1,0} \otimes \Lambda^{2,0}]]$ of real dimension $2\times 4\times 6= 48$ further splits into ${\mathcal W}_1= [[\Lambda^{3,0}]]$ of real dimension $2\times 4=8$ and the remainder ${\mathcal W}_2$. The space $[[\Lambda^{2,1}]]$ splits into another copy of ${\mathcal W}_4=\Lambda^1$, as well as ${\mathcal W}_3=[[\Lambda^{2,1}_0]]$. The last $\Lambda^1$ factor of the intrinsic torsion is denoted by ${\mathcal W}_5$. The first two irreducible factors of the intrinsic torsion ${\mathcal W}_1\oplus {\mathcal W}_2$ control the Nijenhuis tensor, and so vanish when $d\Omega=0$. In addition, $d\Omega=0$ implies the vanishing of the component  ${\mathcal W}_5$ of the intrinsic torsion, which is captured by the $\Lambda^{4,1}$ part of $d\Omega$. 

When $d\Omega=0$, the real part of $d\Phi_\tau=0$ implies $\omega\wedge d\omega=0$. Looking at the decomposition of $\Lambda^5$ into irreducible components, which is formula (3.7) in \cite{Salamon-book}, we see that this equation implies $d\omega=0$. It is known, see \cite{Salamon-book} and also the discussion in \cite{Salamon-Chiossi}, that $d\omega$ captures the components ${\mathcal W}_{1,3,4}$ of the intrinsic torsion. Thus, $d\Omega=0$ together with $d\omega=0$ imply vanishing of all of the intrinsic torsion and thus integrability. 

We have proved the following
\begin{proposition} \label{prop:complex-integrable} The ${\rm SU}(4)$ structure defined by the complex Cayley form $\Phi_\tau$ is integrable if and only if $d\Phi_\tau=0$. 
\end{proposition}

\section{Spinors of ${\rm Spin}(4,4)$ and associated real Cayley forms}
\label{sec:split}

We now describe analogs of all the statements above for the case of the split signature Spin group ${\rm Spin}(4,4)$. This now admits a split octonion model. We deal with corresponding complex Cayley forms in the next section. 

\subsection{Split octonions}

 Split octonions $\Oc'$ form a non-associative normed composition algebra. It is not a division algebra because there are null elements. A split octonion is an object 
\be
\tilde{q} = \tilde{q}_0 \id + \sum_{a=1}^7 \tilde{q}_a {\bf\tilde{e}}^a.
\ee
The unit imaginary octonions ${\bf\tilde{e}}^a$ anti-commute and satisfy
\be\label{split-octonions}
({\bf\tilde{e}}^5)^2=({\bf\tilde{e}}^6)^2=({\bf\tilde{e}}^7)^2=-\id, \quad
({\bf\tilde{e}}^1)^2=({\bf\tilde{e}}^2)^2=({\bf\tilde{e}}^3)^2=({\bf\tilde{e}}^4)^2=\id.
\ee
Thus, the split octonions $\id, {\bf\tilde{e}}^5, {\bf\tilde{e}}^6, {\bf\tilde{e}}^7$ generate a copy of $\Hq\subset\Oc'$. The octonion pairing is given by
\be
(\tilde{q},\tilde{q}) = \tilde{q} \overline{\tilde{q}} = (\tilde{q}_0)^2 - (\tilde{q}_1)^2- (\tilde{q}_2)^2- (\tilde{q}_3)^2- (\tilde{q}_4)^2+(\tilde{q}_5)^2+(\tilde{q}_6)^2+(\tilde{q}_7)^2,
\ee
where the conjugation denoted by overbar changes the signs of all the imaginary generators. 

The products of imaginary octonions are most efficiently encoded into the following 3-form on ${\rm Im}(\Oc')=\R^7$
\begin{equation}\label{tilde-C}
    \tilde{\phi}=e^{567}-e^{5}(e^{41}-e^{23})-e^6(e^{42}-e^{31})-e^7(e^{43}-e^{12}),
\end{equation}
which is similar to (\ref{3-form}) but with the opposite sign in front of the terms containing $e^{1,2,3,4}$. 
This encodes the vector product via 
\be
(u\times v, w) = w \lrcorner v\lrcorner u \lrcorner \tilde{\phi}.
\ee
Here $u\lrcorner$ is the operator of insertion of a vector field $u$ into a differential form. For example ${\bf\tilde{e}}^5 \times {\bf\tilde{e}}^6={\bf\tilde{e}}^7$, but ${\bf\tilde{e}}^5 \times {\bf\tilde{e}}^2=-{\bf\tilde{e}}^3$ because the octonion pairing is negative-definite on directions $1,2,3,4$. 

We note for the later that the 4-form dual to $\phi$ is given by
\be\label{tilde-C-dual}
{}^* \tilde{\phi} = e^{1234} -e^{67}( e^{41}-e^{23}) - e^{75}(e^{42}-e^{31}) - e^{56} (e^{43}-e^{12}).
 \ee

\subsection{The split octonion model of ${\rm Cl}(4,4)$}

Similarly to the case of ${\rm Cl}(8)$, the Clifford algebra ${\rm Cl}(4,4)$ can be realised in terms of $2\times 2$ matrices with entries in split octonions. The general linear combination of Clifford generators is realised as
\be
\Gamma_{\tilde{q}} =\sum_{I=0}^7 \tilde{q}_I \tilde{\Gamma}_I=\left( \begin{array}{cc} 0 & L_{\bar{\tilde{q}}} \\ L_{\tilde{q}} & 0 \end{array}\right),
\ee
where $L_q$ is the operator of left multiplication by $q\in\Oc'$. 

Concretely, an explicit description of the Clifford algebra ${\rm Cl}(4,4)$ in terms of $16\times 16$ matrices arises by taking split octonions to be 8-component columns with $\tilde{q}_0$ as the top entry. One then has the following $\Gamma$-matrices
\begin{equation}\label{Gamma-matr-spin44}
    \tilde{\Gamma}_0=
    \begin{pmatrix}
    0&\id\\
    \id&0
    \end{pmatrix}
    , \quad
    \tilde{\Gamma}_a=
    \begin{pmatrix}
    0&-\tilde{E}_a\\
    \tilde{E}_a&0
    \end{pmatrix}
    \ \text{for} \ a\in \{1,\ldots,7\},
\end{equation}
where
\begin{equation}\label{E-matr}
    \begin{split} 
        \tilde{E}_1&=S_{01}+S_{27}-S_{36}+S_{45},\\
        \tilde{E}_2&=S_{02}-S_{17}+S_{35}+S_{46},\\
        \tilde{E}_3&=S_{03}+S_{16}-S_{25}+S_{47},\\
        \tilde{E}_4&=S_{04}-S_{15}-S_{26}-S_{37}, \\
         \tilde{E}_5&=-E_{05}+E_{14}+E_{23}-E_{67},\\
        \tilde{E}_6&=-E_{06}-E_{13}+E_{24}+E_{57},\\
        \tilde{E}_7&=-E_{07}+E_{12}+E_{34}-E_{56}.
    \end{split}
\end{equation}
Our notation is that $E_{ij}$ is the anti-symmetric matrix that has $+1$ component in the $i$th row and $j$th column, and $-1$ in the $j$th row $i$th column. The matrix $S_{ij}$ is the symmetric matrix with $+1$ in both the $i$th row $j$th column and $j$th row $i$th column positions. The matrices $\tilde{E}_{5,6,7}$ square to minus the identity, while $\tilde{E}_{1,2,3,4}$ square to plus the identity. The matrices $\tilde{E}_a$ anti-commute. 

\subsection{The real split Cayley form}

The canonical split Cayley form is obtained by inserting four distinct $\Gamma$-matrices between two copies of the identity octonion. We get
\be\label{Phi-split}
\tilde{\Phi}:= \langle \id, \tilde{\Gamma}\tilde{\Gamma}\tilde{\Gamma}\tilde{\Gamma}\id \rangle= e^0 \wedge \tilde{\phi} - {}^*\tilde{\phi},
\ee
where $\tilde{\phi}$ is the 3-form (\ref{tilde-C}) that encodes the split octonion multiplication and ${}^*\tilde{\phi}$ is the Hodge dual of $\tilde{\phi}$ given by (\ref{tilde-C-dual}). The stabiliser of $\tilde{\Phi}$ in ${\rm Spin}(3,4)$ is the stabiliser of the Weyl spinor that is the identity split octonion, and this is ${\rm Spin}(3,4)$. 

\subsection{Pure spinors of ${\rm Spin}(4,4)$}

There are several different possible point of view on the real split Cayley form $\tilde{\Phi}$. These arise by representing the real unit spinor (identity split octonion) as a sum of two pure spinors. In contrast to the compact case ${\rm Spin}(8)$, where this decomposition, even though not unique, is possible of only one type, the situation for ${\rm Spin}(4,4)$ is more interesting because there are now pure spinors of different types. 
For Spin groups for indefinite metric signature the structure of the pure spinor orbits is described in \cite{KT}. 

For ${\rm Spin}(4,4)$ there are three different types of pure spinors. They are distinguished by properties of the annihilator subspace. For pure spinors of one type all four vectors spanning $M(\psi_p)$ are complex. In the terminology of reference \cite{KT} these are pure spinors of real index zero. For another type of pure spinors all four vectors are real. These are pure spinors of real index four. And then last, most interesting for our purposes type is that with two real two complex null vectors spanning $M(\psi_p)$. These are pure spinors of real index two. The action of ${\rm Spin}(4,4)$ on all the corresponding orbits is transitive, see \cite{KT}. 

In the previous section, we have seen that a complex viewpoint on the Riemannian real Cayley form is possible, and arises by choosing a complex structure on $\R^8$. At the level of the real unit spinor that goes into the construction of the real Cayley form, this corresponds to writing this unit spinor as a sum of a pure spinor and its complex conjugate. There is a similar viewpoint on the real split Cayley form, except that now there are two different ways that we can represent a real unit split octonion as a sum of two null octonions. 

\subsection{Spit real Cayley form and complex structure}

The first of the two possible viewpoints on the real split Cayley form is the complete analog of what we did in the case of $\R^8$. We choose a complex structure on $\R^{4,4}$, which is encoded by a pure spinor. The corresponding pure spinor $\psi_p$ has four complex null directions spanning $M(\psi_p)$, and is thus of the real index zero. Its complex conjugate $\hat{\psi}_p$ is also a pure spinor, whose null directions $M(\hat{\psi}_p)$ are complementary, i.e. $M(\psi_p)\cap M(\hat{\psi}_p)=0$, and $\langle \psi_p, \hat{\psi}_p\rangle \not=0$. The spinor 
\be\label{psi-real-ind-zero}
\psi = \psi_p+ \hat{\psi}_p
\ee
is not pure, and is a real unit spinor if we normalise $\langle \psi_p, \hat{\psi}_p\rangle =1/2$. Evaluating the 4-form $\Phi=\langle \psi,\Gamma\Gamma\Gamma\Gamma\psi\rangle$ gives us a complex viewpoint on the split real Cayley form. 

Let us make this concrete and consider
\be
\psi_p = \frac{1}{2}(\id + \im \tilde{e}^7),
\ee
which is a complex null split octonion. The corresponding complex structure can be computed by evaluating the 2-form $\langle \overline{\psi}_p, \tilde{\Gamma} \tilde{\Gamma}\psi_p \rangle$ and the 4-form $\langle \psi_p, \tilde{\Gamma} \tilde{\Gamma}\tilde{\Gamma}\tilde{\Gamma}\psi_p \rangle $. This gives
\be
\langle \overline{\psi}_p, \tilde{\Gamma} \tilde{\Gamma}\psi_p \rangle = \frac{\im}{2}\omega, \\ 
\nonumber
\langle \psi_p, \tilde{\Gamma} \tilde{\Gamma}\tilde{\Gamma}\tilde{\Gamma}\psi_p \rangle = \frac{1}{2}\Omega,
\ee
with 
\be\label{omega-Omega-split}
\omega = e^{12}+e^{34}+e^{56}+ e^{70}, \\ \nonumber
\Omega =( e^7+\im e^0)( e^1-\im e^2)( e^3-\im e^4)( e^5+\im e^6).
\ee
The eigenvectors for the corresponding complex structure are
\be
E^+={\rm Span}(e^7+\im e^0, e^1-\im e^2, e^3-\im e^4, e^5+\im e^6 ).
\ee
This is the $+\im$ eigenspace for the operator $R_{\bf{\tilde{e}}^7}$ of right multiplication by the imaginary octonion ${\bf{\tilde{e}}}^7$.

Using $\id = \psi_p+\hat{\psi}_p$ we see that the split Cayley form admits a complex representation
\be\label{Phi-tilde-CS}
\tilde{\Phi} = \langle \id,  \tilde{\Gamma} \tilde{\Gamma}\tilde{\Gamma}\tilde{\Gamma} \id\rangle  = 
\langle \psi_p, \tilde{\Gamma} \tilde{\Gamma}\tilde{\Gamma}\tilde{\Gamma}\psi_p \rangle+ \langle \hat{\psi}_p, \tilde{\Gamma} \tilde{\Gamma}\tilde{\Gamma}\tilde{\Gamma}\hat{\psi}_p \rangle+ 2\langle \hat{\psi}_p, \tilde{\Gamma} \tilde{\Gamma}\tilde{\Gamma}\tilde{\Gamma}\psi_p \rangle=
{\rm Re}({\Omega}) - \frac{1}{2} \omega\wedge \omega,
\ee
with $\Omega,\omega$ given by (\ref{omega-Omega-split}). The original Cayley form is invariant under the group that stabilises the $\id$ spinor, which is ${\rm Spin}(3,4)$. The complex representation makes manifest the ${\rm SU}(2,2)$ subgroup of this stabiliser. 

\subsection{Split real Cayley form and paracomplex structure}

Unlike in the case of ${\rm Spin}(8)$, there is another viewpoint on the real split Cayley form. This arises by considering a pair of complementary real pure spinors $\psi_\pm$. These are now real, and so no longer related by complex conjugation. We will normalise them so that $\langle \psi_+,\psi_-\rangle=1/2$. The stabiliser of such a pair is ${\rm SL}(4,\R)\subset {\rm GL}(4,\R)$, and so such a pair defines an orthogonal paracomplex structure on $\R^{4,4}$, which is a map $K: \R^{4,4}\to \R^{4,4}$ which satisfies $K^2=+\id, g(K\cdot, K\cdot)=-g(\cdot,\cdot)$. The eigenspaces of $K$ are real and totally null. They are identified with the annihilator subspaces $M(\psi_+)$ and $M(\psi_-)$. In addition to a para-complex structure $K$, the pair $\psi_\pm$ also defines top forms on the null eigenspaces of $K$. Writing 
\be\label{psi-real-ind-four}
\psi = \psi_+ + \psi_-
\ee
gives a real unit spinor, and so allows us to get a different perspective on the real split Cayley form $\tilde{\Phi}$, by viewing it from the perspective of a paracomplex structure. 

To make this discussion concrete we choose the following two real pure spinors
\be\label{pure-real}
\psi_\pm = \frac{1}{2}( \id \pm {\bf{\tilde{e}}}^4).
\ee
Each of these pure spinors has the annihilator subspace consisting of four real directions in $\R^{4,4}$. They can be determined by computing $\langle \psi_\pm, \tilde{\Gamma} \tilde{\Gamma}\tilde{\Gamma}\tilde{\Gamma}\psi_\pm \rangle$. We get
\be\label{Omega-psi-pm}
\langle \psi_\pm, \tilde{\Gamma} \tilde{\Gamma}\tilde{\Gamma}\tilde{\Gamma}\psi_\pm \rangle = \frac{1}{2} \Omega_\pm,
\ee
where
\be\label{Omega-pm}
\Omega_+= ( e^4+ e^0)( e^1+ e^5)( e^2+e^6)( e^3+ e^7), \\ \nonumber
\Omega_-= ( e^4- e^0)( e^1- e^5)( e^2-e^6)( e^3- e^7).
\ee
which exhibits the null directions. We also have
\be\label{omega-psi-pm}
\langle \psi_+, \tilde{\Gamma} \tilde{\Gamma}\psi_- \rangle = \frac{1}{2}\omega_r,
\ee
where
\be\label{omega-r}
\omega_r= e^{15} +e^{26}+e^{37}+e^{40}.
\ee
The sum of $\psi_+$ and $\psi_-$ is the original spinor $\id$, and we have the following real representation of the real split Cayley form
\be\label{Phi-split-real}
\tilde{\Phi} = \langle \psi_+, \tilde{\Gamma} \tilde{\Gamma}\tilde{\Gamma}\tilde{\Gamma}\psi_+ \rangle +\langle \psi_-, \tilde{\Gamma} \tilde{\Gamma}\tilde{\Gamma}\tilde{\Gamma}\psi_-\rangle + 2\langle \psi_+, \tilde{\Gamma} \tilde{\Gamma}\tilde{\Gamma}\tilde{\Gamma}\psi_-\rangle =
\frac{1}{2}\Omega_+ +\frac{1}{2} \Omega_- + \frac{1}{2} \omega_r\wedge\omega_r.
\ee
This way of writing the split Cayley form only makes manifest the joint stabiliser of the two pure spinors $\psi_\pm$, which is ${\rm SL}(4,\R)={\rm Spin}(3,3)$.

\subsection{Calibrations in the ${\rm Spin}(4,4)$ case}

There are two different ways of viewing the real split Cayley form (\ref{Phi-split}) from the point of view of calibrations that it defines. Thus, $\tilde{\Phi}$ calibrates subspaces with both Riemannian and split signature metrics. Writing $\tilde{\Phi}$ in a basis adapted to such a calibration exhibits the ${\rm SU}(2)\times{\rm SU}(2)\times{\rm SU}(2)$ and ${\rm SO}(2,1)\times{\rm SO}(2,1)\times{\rm SO}(2,1)$ subgroups of ${\rm Spin}(4,3)$ respectively. 

We first rewrite (\ref{Phi-split}) in a form analogous to (\ref{Phi-calibr}). Thus, we again introduce the self-dual 2-forms in the 4-dimensional subspaces $e^{1,2,3,4}$ and $e^{0,5,6,7}$
\be
\Sigma^1 = e^{41}- e^{23}, \quad \Sigma^2= e^{42}-e^{31}, \quad \Sigma^3= e^{43}-e^{12}, \\ \nonumber
{\Sigma'}^1= e^{05}-e^{67}, \quad {\Sigma'}^1= e^{06}-e^{75}, \quad {\Sigma'}^1= e^{07}-e^{56}.
\ee
The split Cayley form can then be written as 
\be\label{Phi-split-calibr}
\tilde{\Phi} = - \frac{1}{6} \Sigma^i \Sigma^i- \frac{1}{6} {\Sigma'}^i {\Sigma'}^i- \Sigma^i{\Sigma'}^i,
\ee
where the only difference from (\ref{Phi-calibr}) is the sign in front of the last term. This way of writing the Cayley form exhibits the ${\rm SU}(2)\times{\rm SU}(2)\times{\rm SU}(2)$ subgroup of ${\rm Spin}(4,3)$. The subspaces of $\R^{4,4}$ that are locally modelled on $e^{0,5,6,7}$ are calibrations in this case. From the octonionic viewpoint, a choice of such a calibration is a choice of a copy of quaternions $\Hq\subset \Oc'$. 

Another way of writing the same 4-form is to split $\R^8$ into two copies of $\R^{2,2}$. The relevant subspaces in this case are $e^{3,4,5,6}$ and $e^{1,2,7,0}$. We introduce the basis of split signature self-dual 2-forms for these subspaces
\be
\Sigma^1_s = e^{54}+e^{36}, \quad \Sigma^2_s = e^{53}+e^{64}, \quad \Sigma^3_s = e^{56}-e^{43}, \\ \nonumber
{\Sigma'_s}^1 = e^{01} + e^{27}, \quad {\Sigma'_s}^2 = e^{02} + e^{71}, \quad {\Sigma'_s}^3 = e^{07} - e^{12}.
\ee
Then we can write the split Cayley form as
\be\label{Phi-split-calibr-2}
\tilde{\Phi}= 
- \frac{1}{6} \Sigma^i \Sigma^j \eta_{ij}- \frac{1}{6} {\Sigma'}^i {\Sigma'}^j\eta_{ij}+\Sigma^i{\Sigma'}^j \eta_{ij},
\ee
where $\eta_{ij}={\rm diag}(-1,-1,1)$ is the metric on $\R^{2,1}$. This way of writing exhibits the ${\rm SO}(2,1)\times{\rm SO}(2,1)\times{\rm SO}(2,1)$ subgroup of ${\rm Spin}(4,3)$. Thus $\R^{4,4}$ can also be calibrated by copies of subspaces that are locally copies of $\R^{2,2}$. Such a calibration is a choice of a copy of split quaternions $\Hq'\subset \Oc'$.

\section{The Lorentzian Cayley form}
\label{sec:Lorentzian}

We now consider complex Cayley forms for ${\rm Spin}(4,4)$. There are now two one-parameter families of complex Cayley forms.

\subsection{Complex unit spinors for ${\rm Spin}(4,4)$}

Complex Cayley form arises as (\ref{Phi-spinor}) for a complex unit spinor $\psi\in S^+$. Complex unit (Weyl) spinors of ${\rm Spin}(4,4)$ can be classified similarly to the ${\rm Spin}(8)$ case, see \cite{Bhoja:2022txw}.

A complex Weyl spinor is a complexified split octonion $\psi=\alpha_1+\im \alpha_2$. The condition that this spinor is unit becomes $|\alpha_1|^2-|\alpha_2|^2=1, (\alpha_1,\alpha_2)=0$. The pairing between $\psi$ and its complex conjugate $\hat{\psi}\in S^+$ is $\langle \hat{\psi},\psi\rangle= |\alpha_1|^2+|\alpha_2|^2=1+2|\alpha_2|^2$. Because we are dealing with split octonions, this quantity can take any real value. The types of complex spinor orbits are characterised by this value.

When $|\alpha_2|^2>0$ we have two split octonions of positive norm squared. We can then parametrise $|\alpha_1|=\cosh(\tau), |\alpha_2|=\sinh(\tau)$. The spinors
\be
\frac{1}{2} \left( \frac{\alpha_1}{|\alpha_1|} \pm \im \frac{\alpha_2}{|\alpha_2|} \right),
\ee
are both null and thus pure. They are complex conjugates of each other. The stabiliser of either of them is ${\rm SU}(2,2)$, and they define a complex structure on $\R^{4,4}$. Their sum is the real spinor $\alpha_1$, and the more general linear combination (\ref{psi-tau}) gives a complex unit spinor with the same stabiliser. The corresponding complex Cayley form is the precise analog of the complex Cayley form $\Phi_\tau$ in the Riemannian case. 

We note that the orbit with the same stabiliser arises when $|\alpha_2|^2<-1$. In this case both octonions have negative norm squared $|\alpha_1|^2<0, |\alpha_2|^2<0$. One can again introduce two pure spinors that are complex conjugates of each other, so that the complex unit spinor is the linear combination (\ref{psi-tau}). The stabiliser in this case is also ${\rm SU}(2,2)$. 

Another possible orbit is one with $|\alpha_2|^2=0$, in which case $|\alpha_1|^2=1$. Alternatively, one can take $|\alpha_1|^2=0$, and thus $|\alpha_2|^2=-1$. Thus, the complex unit spinor in this case has unit real part, and null imaginary part, or null real part and imaginary part of norm squared minus one. Some discussion on this spinor orbit is contained in \cite{Bhoja:2022txw}, section 5.8. We will not consider this orbit any further. 

The case of most interest to us here is one where one of the two octonions has positive squared norm, while the other has negative norm squared. This is the orbit with $|\alpha_2|^2<0$, while $|\alpha_1|^2>0$, and so $-1<|\alpha_2|^2<0$. The norm squared of $\alpha_1$ then satisfies $0<|\alpha_1|^2<1$. We can parametrise
\be
|\alpha_1|=\cos(\theta), \qquad |\alpha_2|:=\sqrt{-|\alpha_2|^2} = \sin(\theta), \quad \theta\in(0,\pi/2).
\ee
We can then rewrite
\be\label{psi-theta}
\psi_\theta := e^{\im\theta} \psi_+ + e^{-\im\theta} \psi_-,
\ee
where
\be
\psi_\pm := \frac{1}{2}\left( \frac{\alpha_1}{|\alpha_1|} \pm \frac{\alpha_2}{|\alpha_2|}\right)
\ee
are two real null octonions satisfying $\langle \psi_+,\psi_-\rangle=1/2$. We also have $\langle\hat{\psi}_\theta,\psi_\theta\rangle = \cos(2\theta)$. The pair of real null spinors $\psi_\pm$ that satisfy $\langle \psi_+,\psi_-\rangle=1/2$ define a paracomplex structure $K$ on $\R^{4,4}$, together with top forms on the eigenspaces of $K$. The stabiliser of the pair is ${\rm SL}(4,\R)$, which is also the stabiliser of the complex unit spinor $\psi_\theta$. 

\subsection{Complex Cayley form $\Phi_\theta$}

We now build the complex Cayley form using the spinor $\psi_\theta$ (\ref{psi-theta}). All ingredients are already given by the formulas (\ref{Omega-psi-pm}) and (\ref{omega-psi-pm}). We get
\be
\Phi_\theta: = \langle \psi_\theta, \tilde{\Gamma}  \tilde{\Gamma}  \tilde{\Gamma}  \tilde{\Gamma} \psi_\theta \rangle = \cos(2\theta) \frac{1}{2}(\Omega_++\Omega_-) + \im \sin(2\theta) \frac{1}{2}(\Omega_+ - \Omega_-) + \frac{1}{2} \omega_r \wedge \omega_r.
\ee
The stabiliser of this complex Cayley form contains the stabiliser of the pair  of pure spinors $\psi_\pm$, which is ${\rm SL}(4,\R)$. When $\theta=0$ the spinor $\psi$ becomes the sum of two real pure spinors. The corresponding complex Cayley form becomes the real split Cayley form $\Phi_0=\tilde{\Phi}$, whose stabiliser is ${\rm Spin}(4,3)$. When $\theta=\pi/2$ the spinor $\psi_\theta$ is purely imaginary, with the imaginary part given by the difference of two real pure spinors. The corresponding complex Cayley form is again real, with stabiliser ${\rm Spin}(4,3)$. 

\subsection{Lorentzian Cayley form}

When $\theta=\pi/4$ we get the special case $\langle\hat{\psi}_\theta,\psi_\theta\rangle =0$. We call the corresponding complex Cayley form {\bf Lorentzian Cayley form}
\be\label{Lorentzian-Cayley}
\Phi_L := \Phi_{\pi/4} =  \frac{\im}{2}(\Omega_+ - \Omega_-) + \frac{1}{2} \omega_r \wedge \omega_r.
\ee
It is clear that $\Phi_L$ acquires additional discrete symmetry as compared to $\Phi_\theta$, which is to change the sign of both $\Omega_\pm$, followed by the complex conjugation. This is not to be confused with the symmetry of swapping $\Omega_+$ with $\Omega_-$ followed by the complex conjugation, which is a symmetry for any $\theta$, and corresponds to swapping $\psi_\pm$ followed by the complex conjugation. 

Choosing $\psi_\pm$ as in (\ref{pure-real}) we get
 \be\label{psi-compl}
 \psi_L:=\psi_{\pi/4}=\frac{1+\im}{2\sqrt{2}} (\id+{\bf{\tilde{e}}}^4) + \frac{1-\im}{2\sqrt{2}} (\id-{\bf{\tilde{e}}}^4)= \frac{1}{\sqrt{2}}(\id +\im {\bf{\tilde{e}}}^4).
 \ee
We can then write an explicit expression for the resulting Lorentzian Cayley form:
\be
\Phi_L =  \im e^0 \wedge \phi_L + {}^* \phi_L,
\ee
where
\be
\phi_L:= e^{123} - e^1(\im e^{45}-e^{67}) - e^2(\im e^{46}- e^{75})- e^3(\im e^{47}- e^{56}), \\ \nonumber
{}^*\phi_L := \im e^{4567} + e^{23}(\im e^{45}-e^{67}) + e^{31}(\im e^{46}- e^{75})+e^{12}(\im e^{47}- e^{56}).
\ee
Both $\phi_L, {}^*\phi_L$ are built from Lorentzian self-dual 2-forms in directions $e^{4,5,6,7}$, which span a copy of $\R^{1,3}$. The stabiliser of the complex Cayley form is the stabiliser of both $\id, \bf{\tilde{e}}^4$, which is ${\rm Spin}(3,3)={\rm SL}(4,\R)$.

\subsection{A different representation of the Lorentzian Cayley form}

We have seen that the real Riemannian Cayley form is constructed from a unit spinor that can be understood as the sum of two pure spinors, namely a pure spinor and its complex conjugate. Similarly, the real split Cayley form is constructed from a spinor that is given by the sum of a pure spinor $\psi_p$ and its complementary pure spinor $\psi_p'$ such that $\langle \psi_p,\psi_p'\rangle\not=0$. For the real split Cayley form this can be done in two different ways, either as the sum of the pure spinor of real index zero and its complementary, as in (\ref{psi-real-ind-zero}), or as a sum of a pure spinor of real index four and its complementary, as in (\ref{psi-real-ind-four}). 

Thus, the two of the three different types of pure spinors that exist in the split case lead to the split real Cayley form. A natural question to ask is what kind of spinor and the associated Cayley form is obtained when one takes the sum of a pair of pure spinors of the third type. This are pure spinors of real index two, which have two real and two complex null directions. Interestingly, a pure spinor of this type can be understood as a complex linear combination of pure spinors of the other two types. To make this discussion concrete, let us consider the pair
\be\label{pure-pair-22}
\psi_p = \frac{1}{\sqrt{2}}( \frac{1}{2} (\id  + \im {\bf{\tilde{e}}}^7) +\frac{\im}{2} ( {\bf{\tilde{e}}}^4+\im {\bf{\tilde{e}}}^3)) = \frac{1}{\sqrt{2}}( \frac{1}{2} (\id  - {\bf{\tilde{e}}}^3) +\frac{\im}{2} ( {\bf{\tilde{e}}}^4+ {\bf{\tilde{e}}}^7)), \\
\nonumber
\psi_p' = \frac{1}{\sqrt{2}}( \frac{1}{2} (\id  - \im {\bf{\tilde{e}}}^7) +\frac{\im}{2} ( {\bf{\tilde{e}}}^4-\im {\bf{\tilde{e}}}^3)) = \frac{1}{\sqrt{2}}( \frac{1}{2} (\id  + {\bf{\tilde{e}}}^3) +\frac{\im}{2} ( {\bf{\tilde{e}}}^4- {\bf{\tilde{e}}}^7)),
\ee
The way we wrote the two new pure spinors makes it clear that they can be obtained as complex linear combinations of either two pure spinor with four complex null directions, or as a combination of two pure spinors with four real directions. We have $\langle \psi_p',\psi_p\rangle = 1/2$, and so these pure spinors are complementary in that their null subspaces do not intersect. We have
\be
\psi_L= \psi_p+\psi_p',
\ee
and so the complex spinor $\psi_L$ that is used in the construction of the Lorentzian Cayley form can be thought of as obtained as a sum of two complementary pure spinor of real index two. This representation of $\psi_L$ leads to a different formula for $\Phi_L$. 

The null subspace of $\psi_p,\psi_p'$ can be computed by evaluating the corresponding decomposable 4-form. This gives
\be
 \langle \psi_p, \tilde{\Gamma} \tilde{\Gamma}\tilde{\Gamma}\tilde{\Gamma}\psi_p \rangle= \frac{1}{2}\Omega_c, \qquad 
  \langle \psi_p', \tilde{\Gamma} \tilde{\Gamma}\tilde{\Gamma}\tilde{\Gamma}\psi_p' \rangle= \frac{1}{2}\Omega_c',
  \ee
  where
  \be
  \Omega_c = (e^0-\im e^7)(\im e^3+e^4)(e^1+e^6)(e^2-e^5), \\ \nonumber
   \Omega_c' = (e^0+\im e^7)(\im e^3-e^4)(e^1-e^6)(e^2+e^5).
  \ee
  This makes it clear that indeed $\psi_p, \psi_p'$ are of real index two. We can also compute the 2-form
  \be
 \langle \psi_p', \tilde{\Gamma}\tilde{\Gamma}\psi_p \rangle= \frac{1}{2}\omega_c, 
 \ee
 where
 \be
 \omega_c= e^{61}+e^{25}+\im e^{34} -\im e^{07}.
 \ee
 We thus obtain an alternative expression for the Lorentzian Cayley form
 \be\label{Phi-compl-omega}
\Phi_L=  \langle \psi, \tilde{\Gamma} \tilde{\Gamma}\tilde{\Gamma}\tilde{\Gamma}\psi \rangle= 
 \frac{1}{2}(\Omega_c+\Omega_c') + \frac{1}{2} \omega_c\wedge\omega_c.
 \ee

The way of writing (\ref{Phi-compl-omega}) of the complex Cayley form makes manifest the subgroup of ${\rm Spin}(4,4)$ that preserves both $\psi_p, \psi_p'$. This is the subgroup that preserves what the paper \cite{Bhoja:2022dsm} referred to as a structure of a mixed type, which is a complex linear combination of a para-complex structure and the imaginary unit times a complex structure. In this case, this is the group ${\rm SU}(1,1)$ that mixes the complex null directions of the pure spinor $\psi_p$ times ${\rm SL}(2,\R)$ that mixes its real null directions. 

\subsection{Lorentzian calibrations}

We have already seen how the real split Cayley form allows to calibrate $\R^{4,4}$ by 4-dimensional subspaces that are either $\R^4$ or $\R^{2,2}$. We will now exhibit how the complex Cayley form calibrates $\R^{4,4}$ by copies of the Minkowski space. 

We introduce the (complex) basis of self-dual 2-forms in copies of Minkowksi space spanned by $e^{4,5,6,7}$ and $e^{0,1,2,3}$
\be
\Sigma_L^1= \im e^{45}-e^{67}, \quad \Sigma_L^2= \im e^{46}- e^{75}, \quad \Sigma_L^3= \im e^{47}- e^{56}, \\ \nonumber
{\Sigma'_L}^1= \im e^{01}-e^{23}, \quad {\Sigma'_L}^1= \im e^{02}-e^{31}, \quad {\Sigma'_L}^1= \im e^{03}-e^{12}.
\ee
We then have
\be
\Phi_c = -\frac{1}{6} \Sigma^i_L \Sigma_L^i - \frac{1}{6}{\Sigma'_L}^i{\Sigma'_L}^i  - \Sigma_L^i {\Sigma'_L}^i .
\ee
This way of writing the complex Cayley form exhibits the Lorentz subgroup ${\rm Spin}(3,1)$ of its stabiliser ${\rm Spin}(3,3)$. Indeed, this is the diagonal Lorentz subgroup that acts at the same time in both copies of the Minkowski space spanned by $e^{4,5,6,7}$ and $e^{0,1,2,3}$. 

\subsection{Integrable ${\rm SL}(4,\R)$ structures}

With the expression (\ref{Lorentzian-Cayley}) for the Lorentzian Cayley form at hand, it is not difficult to prove that integrable ${\rm SL}(4,\R)$ structures are in one-to-one correspondence with closed Lorentzian Cayley forms.
\begin{proposition} The ${\rm SL}(4,\R)$ structure defined by a Lorentzian Cayley form $\Phi_L$ is integrable if and only if $d\Phi_L=0$. 
\end{proposition}
Indeed, in one direction, an ${\rm SL}(4,\R)$ structure is described by a triple $(g,\omega_r,\Omega^\pm)$, where $g$ is a split signature metric in eight dimensions, the real Kahler form is $\omega_r(\xi_1,\xi_2)=g(K\xi_1,\xi_2)$, where $K$ is an orthogonal para-complex structure $K^2=\id, g(K\xi_1, K\xi_2)=-g(\xi_1,\xi_2)$, and $\Omega^\pm$ are the top forms $\Omega^+\in \Lambda^{4,0}, \Omega^-\in\Lambda^{0,4}$. This structure is integrable if and only if $d\omega_r=0, d\Omega^\pm=0$, which then implies that $\Phi_L$ is closed.

In the other direction, the imaginary part of $d\Phi_L=0$ implies that $d\Omega^\pm=0$, which then implies that the paracomplex structure defined by $\Omega^\pm$ is integrable. The real part of $d\Phi_L=0$ gives $\omega_r d\omega_r=0$. The decomposition of $\Lambda^5 T, T=\R^{4,4}$ into the ${\rm SL}(4,\R)$ irreducible components, formula (3.7) in \cite{Salamon-book}, shows that $\omega_r d\omega_r=0$ implies $d\omega_r=0$, which proves the assertion. We note that the same arguments imply that the more general complex split Cayley forms $\Phi_\tau$ and $\Phi_\theta$ define integrable ${\rm SU}(2,2)$ and ${\rm SL}(4,\R)$ structures respectively if and only if they are closed. 

\section*{Acknowledgements} The author is grateful to N. Hitchin for correspondence, and to Y. Herfary for comments on an early version of this paper.

\section*{Appendix: Lorentzian Cayley forms corresponding to the same metric}

The purpose of this Appendix is to collect some explicit formulas for the tangent space to the space of all Cayley forms corresponding to the same metric. 

\subsection*{Computation of certain 4-forms}

In preparation for the computation of the tangent space to the space of Lorentzian Cayley forms, we first collect some preliminary formulas. 

The base spinor (\ref{psi-compl}) provides us with the operator of the paracomplex structure $K=R_{\tilde{\bf e}^4}$. All quantities are best described by decomposing them into their irreducible components with respect to this para-complex structure. Explicitly, it is given by
\be\label{K}
K = e^4\otimes (e^0)^* + e^0 \otimes (e^4)^* 
+ e^1\otimes (e^5)^* + e^5 \otimes (e^1)^* \\ \nonumber
+e^2\otimes (e^6)^* + e^6 \otimes (e^2)^* +
e^3\otimes (e^7)^* + e^7 \otimes (e^3)^*.
\ee
We note that $\Phi_c$ is invariant under the action of $K$ on all its four indices. 

We now consider the vector space $\R^{4,4}\sim \Oc'$ of real split octonion and decompose it into its components eigenvectors of $K$. Thus, we write
\be
\Xi = \Xi_{1,0} + \Xi_{0,1}, \qquad K(\Xi_{1,0}) = \Xi_{1,0}, \quad K(\Xi_{0,1}) = -\Xi_{0,1}
\ee
The eigenspaces of $K$ are totally null and real. We can parametrise
\be
\Xi_{1,0} =  \xi_{1,0} (\id+{\bf{\tilde{e}}}^4)+ \vec{\xi}_{1,0}, \quad
\Xi_{0,1} =  \xi_{0,1} (\id-{\bf{\tilde{e}}}^4) + \vec{\xi}_{0,1},
\ee
where
\be
\vec{\xi}_{1,0}=  \xi^1_{1,0} ({\bf{\tilde{e}}}^1+{\bf{\tilde{e}}}^5)+ \xi^2_{1,0} ({\bf{\tilde{e}}}^2+{\bf{\tilde{e}}}^6)+ \xi^3_{1,0} ({\bf{\tilde{e}}}^3+{\bf{\tilde{e}}}^6), \\ \nonumber
\vec{\xi}_{0,1}= \xi^1_{0,1} ({\bf{\tilde{e}}}^1-{\bf{\tilde{e}}}^5)+\xi^2_{0,1} ({\bf{\tilde{e}}}^2-{\bf{\tilde{e}}}^6)+ \xi^3_{0,1} ({\bf{\tilde{e}}}^3-{\bf{\tilde{e}}}^6).
\ee
For purposes of the next subsection, we list the following results, obtained by explicit computation:
\be\label{res}
\langle \psi_+, \tilde{\Gamma} \tilde{\Gamma} \tilde{\Gamma} \tilde{\Gamma} \Xi_{1,0}\rangle = \Omega_+ \xi_{1,0} + \omega_r \wedge \vec{\xi}_{1,0} \lrcorner \omega_+, \\ \nonumber
\langle \psi_+, \tilde{\Gamma} \tilde{\Gamma} \tilde{\Gamma} \tilde{\Gamma} \Xi_{0,1}\rangle = \frac{1}{2} \omega_r \omega_r \xi_{0,1} 
+ \omega_r \wedge (e^4+e^0) \wedge \vec{\xi}_{0,1} \lrcorner \omega, 
\\ \nonumber
\langle \psi_-, \tilde{\Gamma} \tilde{\Gamma} \tilde{\Gamma} \tilde{\Gamma} \Xi_{1,0}\rangle = \frac{1}{2} \omega_r \omega_r \xi_{1,0} 
+ \omega_r\wedge  (e^4-e^0) \wedge \vec{\xi}_{1,0}\lrcorner \omega,
\\ \nonumber
\langle \psi_-, \tilde{\Gamma} \tilde{\Gamma} \tilde{\Gamma} \tilde{\Gamma} \Xi_{0,1}\rangle = \Omega_- \xi_{0,1}  + \omega_r\wedge \vec{\xi}_{0,1}\lrcorner \omega_-,
\ee
where we introduced the 1- and 2-forms
\be
\vec{\xi}_{0,1} \lrcorner \omega:= (e^1+e^5) \xi^1_{0,1} + (e^2+e^6) \xi^2_{0,1} + (e^3+e^7) \xi^3_{0,1}, \\ \nonumber
\vec{\xi}_{1,0} \lrcorner \omega_+ :=  (e^2+e^6)(e^3+e^7) \xi^1_{1,0} + (e^3+e^7) (e^1+e^5)\xi^2_{1,0} + (e^1+e^5)(e^2+e^6) \xi^3_{1,0}, \\ \nonumber
\vec{\xi}_{1,0}\lrcorner \omega := (e^1-e^5) \xi^1_{1,0} + (e^2-e^6) \xi^2_{1,0} + (e^3-e^7) \xi^3_{1,0}, \\ \nonumber
\vec{\xi}_{0,1}\lrcorner \omega_-:= (e^2-e^6)(e^3-e^7) \xi^1_{0,1} + (e^3-e^7) (e^1-e^5)\xi^2_{0,1} + (e^1-e^5)(e^2-e^6) \xi^3_{0,1}.
\ee
Here we are referring to the following 2- and 3-forms:
\begin{align}
&\omega:= e^{15}+e^{26}+e^{37}, \\ \nonumber
&\omega_+ = (e^1+e^5) (e^2+e^6) (e^3+e^7) , \qquad \omega_- = (e^1-e^5) (e^2-e^6) (e^3-e^7) .
\end{align}

\subsection*{Lorentzian Cayley forms corresponding to the same metric}

In the positive definite case, the space of Cayley forms giving rise to the same metric is the 7-dimensional space of unit spinors. We have seen that such Cayley forms can be parametrised explicitly as (\ref{fibre-unit-spinor}). We now want to obtain a similar parametrisation in the case of complex Cayley forms. 

The first step is to parametrise the space of complex unit spinors that are orthogonal to their complex conjugates. Such a parametrisation is provided by (\ref{psi-theta}) with $\theta=\pi/4$. The spinor $\psi=\psi_{\pi/4}$ arising this way is automatically orthogonal to its complex conjugate provided both $\psi_\pm$ are real and pure (and thus null). The norm of $\psi$ is then $\langle\psi,\psi\rangle= 2\langle \psi_+,\psi_-\rangle$, which is automatically real. We get unit spinors $\psi$ if we impose $\langle \psi_+,\psi_-\rangle=1/2$.

We now compute the tangent space to the set of such spinors, at $\psi_\pm$ given by (\ref{pure-real}). The condition that each $\psi_\pm$ remains null becomes
\be
\langle \psi_+, \delta \psi_+\rangle =0, \qquad \langle \psi_-, \delta \psi_-\rangle =0, 
\ee
which means that 
\be\label{psi-pm-pert}
\delta \psi_+ = a (\id+{\bf{\tilde{e}}}^4) + \vec{\xi}, \qquad
\delta \psi_- = b (\id-{\bf{\tilde{e}}}^4) + \vec{\eta},
\ee
where $a,b\in \R$ and  $\vec{\xi},\vec{\eta}\in {\rm Span}( {\bf{\tilde{e}}}^{1,2,3,5,6,7})$. The condition that $\langle \psi_+,\psi_-\rangle=1/2$ remains satisfied implies
\be
\langle \psi_+, \delta \psi_-\rangle + \langle \delta \psi_+, \psi_-\rangle=0,
\ee
which implies $a+b=0$ in (\ref{psi-pm-pert}). This gives the expected dimension $6+6+1$ for the tangent space. 

We now take 
\be
\delta \psi = \frac{1+\im}{\sqrt{2}} \delta\psi_+ + \frac{1-\im}{\sqrt{2}} \delta\psi_-,
\ee
and compute
\be\label{delta-phi-c}
2\delta\Phi_c = \langle \psi, \tilde{\Gamma} \tilde{\Gamma}\tilde{\Gamma}\tilde{\Gamma} \delta\psi\rangle =
\im \langle \psi_+, \tilde{\Gamma} \tilde{\Gamma}\tilde{\Gamma}\tilde{\Gamma} \delta\psi_+\rangle
-\im \langle \psi_-, \tilde{\Gamma} \tilde{\Gamma}\tilde{\Gamma}\tilde{\Gamma} \delta\psi_-\rangle
\\ \nonumber
+ \langle \psi_+, \tilde{\Gamma} \tilde{\Gamma}\tilde{\Gamma}\tilde{\Gamma} \delta\psi_-\rangle
+\langle \psi_-, \tilde{\Gamma} \tilde{\Gamma}\tilde{\Gamma}\tilde{\Gamma} \delta\psi_+\rangle.
\ee

We first compute the part of the tangent vector in the only direction that involves changing the $\id, {\bf{\tilde{e}}}^4$ components of the spinor. Thus, we consider
\be
\delta \psi_+ = a (\id+{\bf{\tilde{e}}}^4), \qquad
\delta \psi_- = -a (\id-{\bf{\tilde{e}}}^4).
\ee
We use the results (\ref{res}). For this tangent vector the real parts in (\ref{delta-phi-c}) cancel and we have
\be\label{delta-Phi-1}
2\delta \Phi_c = \im a( \Omega_+ + \Omega_-).
\ee
We then compute the variation of the Cayley form in the remaining tangent directions spanned by $\vec{\xi},\vec{\eta}$. We get
\be\label{delta-Phi-2}
2\delta \Phi_c = \im \omega_r \wedge\left( (  \vec{\xi}_{1,0} \lrcorner \omega_+ +  (e^4+e^0)\wedge \vec{\xi}_{0,1} \lrcorner \omega) 
- (  \vec{\eta}_{0,1} \lrcorner \omega_- +  (e^4-e^0)\wedge \vec{\eta}_{1,0} \lrcorner \omega)\right)\\ \nonumber
+\omega_r \wedge \left( (\vec{\eta}_{1,0} \lrcorner \omega_+ +  (e^4+e^0)\wedge \vec{\eta}_{0,1} \lrcorner \omega) 
+(  \vec{\xi}_{0,1} \lrcorner \omega_- + (e^4-e^0)\wedge \vec{\xi}_{1,0} \lrcorner \omega)\right).
\ee
Here the first two terms in brackets in the first line span all of $\Lambda^{2,0}$, and the second two terms span all of $\Lambda^{0,2}$. The same spaces are spanned in the brackets in the second line, but differently parametrised. 

We can collect these results as the following proposition. Let $\Omega^4_c(M)\subset\Lambda^4_\C(M)$ be the space of Lorentzian Cayley forms on an 8-manifold $M$. Every such form $\Phi_c$ defines a split signature metric. Let us fix $\Phi_c$, and let $g=g_{\Phi_c}$ be the (real) metric defined by $\Phi_c$. Let $\Sigma_g$ be the bundle of all Lorentzian Cayley forms that give rise to the same metric. This is a real rank $13$ bundle, that, as we know, can be thought of as the bundle of complex unit spinors that are orthogonal to their complex conjugates. Alternatively, each fibre of $\Sigma_g\to M$ is a copy of ${\rm Spin}(4,4)/{\rm Spin}(3,3)$. Then the original Lorentzian Cayley form $\Phi_c$ is a section of $\Sigma_g\to M$. 

\begin{proposition} Let $S_{\Phi_c}$ be the tangent space to the space of Cayley forms corresponding to the same metric, i.e. the vector bundle, subbundle of $\Lambda^4_\C(M)$, whose fibre at each point $x\in M$ is the tangent to the fibre $\Sigma_g$ at $\Phi_c$. Then $S_{\Phi_c}$ decomposes into the following irreducible components
\be\label{tangent-space-decomp}
{\rm Re}(S_{\Phi_c}) = \Lambda^{2,0} \oplus \Lambda^{0,2}, \\ \nonumber
{\rm Im}(S_{\Phi_c}) = \Lambda^{2,0} \oplus \Lambda^{0,2} \oplus \R,
\ee
with a bijective map between the spaces $\Lambda^{2,0}, \Lambda^{0,2}$ appearing in the real and imaginary parts of $S_{\Phi_c}$. Explicitly, the relevant factors are described by the formulas (\ref{delta-Phi-1}) and (\ref{delta-Phi-2}).
\end{proposition}

We remark that the spaces that appear in the decomposition of $S_{\Phi_c}$ are as expected from the general theory of intrinsic torsion. Indeed, given an ${\rm SL}(4,\R)$ structure, the tangent space $T=\R^{4,4}$ splits as $\Lambda^{1,0}\oplus \Lambda^{0,1}$. The Lie algebra $\mathfrak{spin}(8)=\Lambda^2 T$ then splits as
\be
\mathfrak{spin}(8)= \Lambda^{2,0}\oplus \Lambda^{0,2} \oplus \Lambda^{1,1}_0 \oplus \R.
\ee
Here the $\Lambda^{1,1}_0$ factor is 
\be
\mathfrak{sl}(4,\R) = \Lambda^{1,1}_0,
\ee
and so
\be
\mathfrak{sl}(4,\R)^\perp = \Lambda^{2,0}\oplus \Lambda^{0,2} \oplus \R.
\ee
These are precisely the spaces that appear in (\ref{tangent-space-decomp}). The above proposition describes how these spaces are embedded into the space $\Lambda^4_\C(\R^{4,4})$ of complexified 4-forms.



\end{document}